\documentclass[11pt]{article}

\usepackage{amsfonts}
\usepackage{amssymb}
\usepackage{graphicx}

\usepackage{ifthen}
\usepackage{color}
\usepackage{cite}
\input{steve.lib}

\begin{document}

\title{Exponential Convergence of a
Distributed Algorithm for Solving Linear Algebraic Equations\thanks{
The material in this paper was partially presented at the 53rd IEEE Conference on Decision and Control \cite{cdc14.2}.
The authors wish to thank Anup Rao (Georgia Institute of Technology)
and Meiyue Shao (Lawrence Berkeley National Laboratory)
for useful discussions which have contributed to this work.
J. Liu is with University of Illinois at Urbana-Champaign (\texttt{jiliu@illinois.edu}).
A. S. Morse is with Yale University, USA (\texttt{as.morse@yale.edu}).
A. Nedi\'c is with Arizona State University (\texttt{angelia.nedich@asu.edu}).
T. Ba\c{s}ar is with University of Illinois at Urbana-Champaign (\texttt{basar1@illinois.edu}).
}}

\author{Ji Liu  \hspace{.2in}  A. Stephen Morse  \hspace{.2in}  Angelia Nedi\'c  \hspace{.2in}  Tamer Ba\c{s}ar}
\maketitle

\begin{abstract}
In a recent paper, a distributed algorithm was proposed for solving linear algebraic
equations of the form $Ax = b$ assuming that the equation has at least one solution. The
equation is presumed to be solved by $m$ agents assuming that each agent knows a subset of the
rows of the matrix $\matt{A & b}$, the current estimates of the equation's solution generated by each
of its neighbors, and nothing more. Neighbor relationships are represented by a time-dependent
directed graph $\bbb{N}(t)$ whose vertices correspond to agents and whose arcs characterize neighbor
relationships. Sufficient conditions on $\bbb{N}(t)$ were derived under which the algorithm can cause all
agents' estimates to converge exponentially fast to the same solution to $Ax = b$. These
conditions were also shown to be necessary for exponential convergence, provided the data
about $\matt{A & b}$ available to the agents is ``non-redundant''. The aim of this paper is to relax
this ``non-redundant'' assumption. This is accomplished by establishing exponential convergence
under conditions which are the weakest possible for the problem at hand;
{\color{black} the conditions are based on a new notion of graph connectivity.
An improved bound on the convergence rate is also derived.}
\end{abstract}


\section{Introduction}

{\color{black}
Over the past few decades, there has been considerable interest in developing algorithms for
information distribution and computation among agents via local interactions
\cite{borkar,Ts3,Ts4,ber,lynch,vicsekmodel}.
Recently, the need for distributed processing has arisen naturally in multi-agent
and sensor networks \cite{survey,moura,pieee} because autonomous agents or mobile sensors are physically separated from each other
and communication constraints limit the flow of information across a multi-agent or sensor network
and consequently preclude centralized processing.
As a consequence, distributed computation and decision making problems of all types have arisen naturally;
notable examples include consensus, multi-agent coverage problems, the rendezvous problem,
localization of sensors in a multi-sensor network, and the distributed management of multi-agent
formations.
One of the most important numerical computations involving
real numbers is solving a system of linear algebraic equations,
which has received much attention for a long time, especially in the parallel
processing community where the main objective is to solve the system faster or more accurately.
It is with these thoughts in mind that we are interested in the problem of solving a system of linear algebraic equations in a distributed manner,
introduced in more precise terms as follows.
}

Consider a network of $m>1$ autonomous agents
 which are able to  receive information from  their ``neighbors''.
  Neighbor relationships  are
  characterized by a time-dependent directed graph $\bbb{N}(t)$ with $m$ vertices and a set of
   arcs defined so that    there is an arc in the graph
   from vertex $j$ to vertex $i$ whenever agent $j$ is a {\em neighbor} of agent $i$. Thus,  the  directions
    of arcs represent the directions of
   information flow.
For simplicity, we take each agent to be a neighbor of itself.
Thus, $\bbb{N}(t)$ has self-arcs at all vertices.
Each agent $i$ has a real-time  dependent state vector $x_i(t)$ taking values in $\R^n$,
 and we assume that the information agent
 $i$ receives from neighbor $j$ is only the current state vector of neighbor $j$.
  We also assume that agent $i$
 knows only a pair of real-valued  matrices $(A_i^{n_i\times n},b_i^{n_i\times 1})$.
The problem of interest is to devise local algorithms,
 one for each agent, which will enable all $m$ agents to iteratively compute
   the same solution
 to the linear equation $Ax=b$
where
$$A = \matt{A_1 \cr A_2 \cr \vdots \cr A_m}_{\bar{n}\times n},  \;\;\;\;\;\;\;\;
b= \matt{b_1\cr b_2\cr \vdots \cr b_m}_{\bar{n}\times 1}$$
and  $\bar{n} =\sum_{i=1}^m n_i$.
{\color{black}
We assume that  $Ax=b$  has  at least one solution,
unless stated otherwise.
The algorithm presented in this paper works for both the case when $Ax=b$
has a unique solution and the case when $Ax=b$ has multiple solutions.
For the case when $Ax=b$ does not have a solution, the algorithm can be modified
to obtain a least squares solution via a centralized initialization step
(see Section \ref{least}).
}


Recently, a distributed algorithm was proposed in \cite{cdc13.1}
for the  synchronous version of the problem just formulated, and with slight modification,
that is for a restricted asynchronous version of the problem in which transmission delays are not taken into account.
A more general asynchronous version of the problem in which
transmission delays are explicitly taken into account was addressed in \cite{cdc13}.

{\color{black}
The synchronous version of the problem considered here can be viewed as a distributed parameter estimation problem
\cite{xxx,infotheory,calibration}.
One approach to the problem is to reformulate it as a distributed convex optimization problem,
which has a rich literature \cite{nedic,nedic2,nedic4,xavier,wainwright,cortes,wei,nedic5,review12,review14,reviewlong}.
An alternative approach to the problem 
is to view it as a constrained consensus problem \cite{nedic2,ren,ren14,cdc14.3}.
A similar problem with more restrictive assumptions has been studied in \cite{lu3,lu1}.
The problem is related to classical parallel algorithms
such as Jacobi iterations \cite{margaris}, so-called ``successive over-relaxations''
   \cite{sor}, and the Kaczmart method \cite{kaczmarz}.
The problem is also related to the problem of estimation on graphs from relative measurements
in which $A$ is determined by the underlying graph and noisy measurements are taken into account
\cite{estimation07,estimation08,estimation09}.
}

The differences and advantages of the algorithm in \cite{cdc13.1}, compared with those
in the literature \cite{nedic2,ren,xxx,calibration,tron,nedic}
\cite{cortes,lu3,lu1,infotheory,margaris,sor,kaczmarz},
have been discussed in \cite{cdc13.1,cdc13,tac}.
{\color{black}
Specifically,
the algorithm in \cite{cdc13.1}
\begin{enumerate}
  \item is applicable
to any pair of real matrices $(A, b)$ for which $Ax = b$ has at
least one solution,
  \item is capable of finding a solution
exponentially fast,
  \item is capable of finding a solution for a time-varying directed graph sequence
under appropriate joint connectedness,
  \item is capable of finding
a solution using at most an $n$-dimensional state vector received
at each clock time from each of its neighbors,
  \item is applicable
without imposing restrictive requirements such
as a) the assumption that each agent is constantly aware of
an upper bound on the number of neighbors of each of its
neighbors or b) the assumption that all agents are able to share
the same time-varying step size.
\end{enumerate}
See Section II in \cite{tac} for details.
To the best of our knowledge, there is no distributed convex optimization algorithm
which simultaneously satisfies all the above properties.
We provide a comparison with competing algorithms in the following table.
\begin{center}
    \begin{tabular}{ | p{1.8cm} |  p{2.8cm} |  p{2.5cm} |}
    \hline
    Paper & Convergence Rate & Neighbor Graph \\ [5pt]\hline
    this paper & exponentially fast & time-varying, directed \\ [5pt]\hline
    \cite{calibration} & exponentially fast & time-varying, undirected \\ [5pt]\hline
    \cite{nedic2} & exponentially fast & time-invariant, complete \\ [5pt]\hline
    \cite{nedic5} & $O(\ln t / \sqrt t)$ & time-varying, directed \\ [5pt]\hline
    \cite{review15cdc} & (locally) exponentially fast & time-invariant, directed \\ [5pt]\hline
    \cite{review16tsp} & exponentially fast & time-invariant, star \\ [5pt]\hline
    \cite{review14hk} & $O(1/t)$ & time-invariant, star \\ [5pt]\hline
    \cite{review12} & not explicit & time-invariant, undirected \\ [5pt]\hline
    \end{tabular}
\end{center}
From the table, it can be seen that only the algorithm presented in this paper
can solve the problem exponentially fast for time-varying, directed, neighbor graphs.
It is worth noting that the idea in \cite{calibration} can
solve the problem for time-varying, directed, neighbor graphs by using double linear iterations
which are specifically tailored to the distributed averaging problem when unidirectional
communications exist \cite{acc12};
but the downside of this idea is that the
amount of data to be communicated between agents does not
scale well as the number of agents increases.
}

{\color{black}
Continuous-time distributed algorithms for the problem in this paper have also received some attention lately; see \cite{tang,brian,shi,acc15}.
}

{\color{black}
From the preceding discussion, a significant advantage of
the algorithm in \cite{cdc13.1} over the other existing ones is its capability to
solve the problem exponentially fast
even when the underlying neighbor graph is directed and time-varying,
using only an $n$-dimensional state vector transmitted between neighboring agents at each clock time.
Accordingly, our aim in this paper is to analyze the algorithm proposed in \cite{cdc13.1},
and particularly to determine the weakest graph-theoretic condition under which
the algorithm can solve the distributed linear equation problem
exponentially fast.
We emphasize exponential convergence because it is robust against certain types of perturbation,
analogous to exponential stability of linear systems \cite{rugh}; it will be clear shortly that
the system determined by the algorithm in \cite{cdc13.1} is a discrete-time linear time-varying system.
}

In this paper, we focus on the synchronous version of the problem,
but the results derived can be straightforwardly extended to asynchronous versions using the analysis tools in \cite{cdc13}.
In \cite{tac}, a necessary and sufficient graph-theoretic condition was obtained under a ``non-redundant'' assumption.
{\color{black}
Roughly speaking, the set of $m$ agents is non-redundant if a distributed solution to $Ax = b$ cannot
be obtained by any proper subset of the full set of $m$ agents; otherwise, the set is redundant.
The formal definition is given as follows.

We say
that agents with  labels in
$\scr{V} =\{i_1,i_2,\ldots, i_q\}\subset \{1,2,\ldots,m\}$  are {\em redundant} if any solution to the equations
$A_ix = b_i$ for  all $i$ in the complement of $\scr{V}$, is a solution to $Ax = b$.
To derive an algebraic condition for redundancy,
suppose that $z$ is a solution to $Ax=b$. Write $\bar{\scr{V}}$
 for the complement of $\scr{V}$ in $\{1,2,\ldots,m\}$.   Then, any solution $w$ to the equations $A_ix=b_i$, $i\in\bar{\scr{V}}$,
 must satisfy $w-z\in\bigcap_{i\in\bar{\scr{V}}} \ {\rm ker} \ A_i$.
 Thus, agents with labels in  $\scr{V}$ will be redundant
  whenever  $w-z\in\bigcap_{i\in\scr{V}}\ {\rm ker} \ A_i$.
Therefore, agents with labels in  $\scr{V}$ will be redundant if, and only if,
$$\bigcap_{i\in\bar{\scr{V}}}\ {\rm ker} \ A_i\subset \bigcap_{i\in\scr{V}}\ {\rm ker} \ A_i$$
We say that agents with labels in $\{1,2,\ldots, m\}$ is a {\em non-redundant} set if no such proper subset exists.
}

{\color{black}
Suppose that $A\neq 0$ and the data about $\matt{A & b}$ available to the agents is non-redundant;
then it has been shown in \cite{tac} that the synchronous algorithm causes all agents' estimates to converge exponentially fast
to the same solution to $Ax=b$ if, and only if, the sequence of neighbor graphs is
``repeatedly jointly strongly connected''.
}
{\color{black}
Since the agents acquire the data about $\matt{A & b}$ in a distributed manner,
without coordination among the agents,
it cannot be guaranteed that the data available to the agents is non-redundant.
Thus,}
the following questions remain.
Is it possible to relax the non-redundancy assumption?
What is the weakest possible graph-theoretic condition for exponential convergence under more general assumptions?
It is with these questions in mind that we investigate in depth the stability of the algorithm.

The main contribution of this paper is to provide necessary and sufficient graph-theoretic conditions
for the algorithm to converge exponentially fast without the non-redundancy assumption,
{\color{black}
which are the weakest possible conditions for exponential convergence.
}
The conditions are based on a parameter-dependent notion of graph connectivity, which is
less restrictive than strong connectivity,
and thus generalize the results in \cite{tac}.
{\color{black}
A convergability issue, a further result on the new notion of connectivity, an improved
bound on convergence rate, and least squares solutions, are also addressed and discussed (in Section \ref{discuss}).
}

The material in this paper was partially presented in \cite{cdc14.2},
but this paper presents a more comprehensive treatment of the work.
Specifically, the paper provides proofs for Propositions \ref{ns} and \ref{p}, Lemmas \ref{eig1} and \ref{vlad}
(Lemma 2 in \cite{cdc14.2}),
and
establishes additional results in \S \ref{discuss}, which were not included in \cite{cdc14.2}.

\subsection{Organization}

The remainder of this paper is organized as follows.
Some preliminaries are provided in \S \ref{notion}.
The synchronous algorithm proposed in \cite{cdc13.1} is revisited in \S \ref{system}.
All $m$ agent update rules can be combined into one linear time-variant system
whose update matrices are related to ``flocking matrices'' in consensus problems \cite{vicsekmodel},
but are of a more complicated form.
Thus, the problem can be viewed as a generalized but more challenging consensus problem.
We first treat in \S \ref{pconnect} a special case when the linear equation $Ax=b$ has a unique solution
and the neighbor graph $\bbb{N}(t)$ is independent of time,
which simplifies to a linear time-invariant system;
after establishing the necessary and sufficient graph-theoretic condition for the system to be
exponentially stable, we introduce a new notion of graph connectivity.
Based on this notion, we state the main results of this paper in \S \ref{main}.
In particular, necessary and sufficient graph-theoretic conditions are presented for
both the unique solution case (Theorem \ref{main1}) and the nonunique solution case (Theorem \ref{main2}),
which provide the weakest possible conditions for the algorithm to converge exponentially fast.
The proofs of the main results are given in \S \ref{analysis}.
The uniqueness case is analyzed first in \S \ref{unique} and the nonuniqueness case is treated next in \S \ref{nonunique}.
In the latter case, we appeal to a result on the stability of linear consensus processes \cite{luc,cdc14}.
Finally, some additional results are discussed in \S \ref{discuss}.
A convergability issue is addressed in \S \ref{compare},
a further result on the new notion of graph connectivity is given in \S \ref{equal},
{\color{black}
an improved bound on convergence rate is derived in \S \ref{rate},
and the case when $Ax=b$ does not have a solution is discussed in \S \ref{least}
for obtaining a least squares solution.
}

\subsection{Preliminaries} \label{notion}

If $n$ is a positive integer, we define $\mathbf{n}=\{1,2,\ldots,n\}$.
For a set of $m$ matrices $\{M_1,M_2,\ldots,M_m\}$ with the same number of columns, we define
$${\rm stack} \{M_1,M_2,\ldots,M_m\}= \matt{M_1 \cr M_2 \cr \vdots \cr M_m}$$
A nonnegative $m\times m$ matrix is called {\em stochastic} if its row sums all equal $1$.
We use $\scr{G}_{sa}$ to denote the set of all directed graphs with
  $m$ vertices 
    which have self-arcs at all vertices.
The graph of an $m\times m$ nonnegative matrix  $M$ is an $m$-vertex directed
 graph $\gamma(M)$  defined so that $(i,j)$ is an arc from $i$ to $j$  in the graph whenever
  the $ji$th entry of
 $M$ is nonzero. Such a graph will be in $\scr{G}_{sa}$ if and only if all diagonal entries of
  $M$ are positive.

Let $\bbb{G}_p$ and $\bbb{G}_q$ be two directed graphs with $m$ vertices.
  By the {\em composition} of $\mathbb{G}_p$
with $\mathbb{G}_q$,  denoted by $\mathbb{G}_q\circ\mathbb{G}_p$,
is meant that directed graph with $m$ vertices and
 arc set defined  so that $(i, j)$ is an arc in the
composition whenever there is a vertex $k$ such that $(i, k)$
 is an arc in $\mathbb{G}_p$ and $(k, j)$
is an arc in $\mathbb{G}_q$.
Note that composition is an associative binary operation; because of this,
the definition extends unambiguously to any finite sequence of directed graphs with the same vertex set.
Composition is defined so that
for any pair of  $m\times m$ nonnegative matrices $M_1$ and $M_2$, there holds
 $\gamma(M_2M_1)=\gamma(M_2)\circ\gamma(M_1)$.
If we focus exclusively on graphs in $\scr{G}_{sa}$, more can be said. In this case,
the definition implies that the arcs of both $\bbb{G}_p$ and $\bbb{G}_q$ are arcs of $\bbb{G}_q\circ\bbb{G}_p$;
the converse is false.

A directed graph $\mathbb{G}$ is {\em strongly connected} if there is a directed path between
each pair of distinct vertices. A directed graph $\mathbb{G}$ is {\em rooted} if
it contains a directed spanning tree of $\bbb{G}$. Note that every strongly connected graph is rooted;
the converse statement is false.
We say that a finite sequence of directed graphs $\bbb{G}_1,\bbb{G}_2,\ldots,\bbb{G}_p$ with the same vertex set
is {\em jointly strongly connected} (or {\em jointly rooted}) if the composition
$\mathbb{G}_{q}\circ\mathbb{G}_{q-1}\circ  \cdots\circ\mathbb{G}_{1}$ is strongly connected (or rooted).
We say that an infinite  sequence of directed graphs $\mathbb{G}_1, \mathbb{G}_2,\ldots$ with the same vertex set
is {\em repeatedly jointly strongly connected} (or {\em repeatedly jointly rooted})
if there exist finite positive integers $l$ and $\tau_0$ such that for any integer $k\geq 0$,
the finite sequence  $\mathbb{G}_{\tau_0+kl},\mathbb{G}_{\tau_0+kl+1},\ldots,\mathbb{G}_{\tau_0+(k+1)l-1}$
is jointly strongly connected (or jointly rooted).
These notions of connectivity are more or less well known in the study of
distributed averaging and consensus problems \cite{acc12,reachingp1},
although the form of the condition may vary slightly from publication to publication.
See for example \cite{nedic3,luc}.

\section{The Algorithm}\label{system}

The system to be studied consists of $m>1$ autonomous agents labeled $1$ through $m$.\footnote{
The purpose of labeling is merely for convenience in expressions. We do
not require such a global ordering.}
As mentioned in the introduction, we are interested in the problem of solving linear equations of the form
$Ax=b$, where $A$ is a matrix and $b$ is a vector, for which the equation has at least one solution,
in a distributed manner among the $m$ agents.
Suppose that time is discrete in that $t$ takes values in the set $\{1,2,\ldots\}$.
The synchronous algorithm proposed in \cite{cdc13.1} is as follows.
Each agent $i$ initializes its state $x_i\in\R^n$ at time $t=1$ by picking $x_i(1)$ to be any solution
to the equation $A_ix=b_i$. From then on, each agent $i$
iteratively updates its state using
$$x_i(t+1) = x_i(t)- \frac{1}{m_i(t)}P_i\left(m_i(t)x_i(t)-\sum_{j\in\scr{N}_i(t)}x_j(t)\right),$$
\eq{i\in\mathbf{m}, \;\;\;\;\; t\geq 1\label{a1}}
where $\scr{N}_i(t)$ is the set of labels of agent $i$'s neighbors at time $t$,
$m_i(t)$ is the number of labels in $\scr{N}_i(t)$ (or equivalently, the in-degree of vertex $i$ in $\bbb{N}(t)$),
and $P_i$ is the readily computable orthogonal projection on the kernel of  $A_i$.
{\color{black}
The key idea behind the algorithm \rep{a1} is as follows.
  Suppose that
$K_i$ is a basis matrix for the kernel of $A_i$.
If we restrict the updating  of $x_i(t)$ to iterations of the form
$x_i(t+1) = x_i(t) + K_iu_i(t)$, $t\geq 1$,
then no matter what $u_i(t)$ is,  each $x_i(t)$ will satisfy $A_ix_i(t) = b_i$, $t\geq 1$.
Then, in accordance with
the agreement principle, all we need to do to solve the problem
is to come up with a good way to choose the $u_i$ so that a consensus
is ultimately reached.
The idea here is to choose $x_i(t+1)$
to satisfy $A_ix_i(t+1) = b_i$  while at the same time making  $x_i(t+1)$  approximately
 equal to the average of agent $i$'s neighbors'  current estimates of the solution to $Ax=b$.
Specifically, we choose each $u_i(t)$ to minimize the difference
 $\left (x_i(t)+K_iu_i(t)\right ) - \frac{1}{m_i(t)}\left(\sum_{j\in\scr{N}_i(t)}x_j(t)\right )$ in the least squares sense.
Doing this leads at once to an iteration for
 agent $i$
  of the form \rep{a1}.
See Section III of \cite{cdc13.1} for more details.
}

It is possible to combine the above $m$ update equations into one linear recursion equation.
Toward this end, let $x^*$ be a solution to $Ax=b$ and define
$$y_i(t) = x_i(t)-x^*,\;\;\;\;\; i\in\mathbf{m},\;\;\;\;\; t\geq 1$$
Then, it has been shown in \cite{cdc13.1} that \rep{a1} simplifies to
\eq{y_i(t+1) = \frac{1}{m_i(t)}P_i\sum_{j\in\scr{N}_i(t)}P_jy_j(t),\;\;\; i\in\mathbf{m},\;\;\; t\geq 1\label{na1}}
Set $y(t) = {\rm stack}\{y_1(t),y_2(t),\ldots, y_m(t)\}$.
Let
 $A_{\mathbb{N}(t)}$ denote the adjacency matrix of $\mathbb{N}(t)$,  $D_{\mathbb{N}(t)}$
 denote the $m\times m$ diagonal matrix whose $i$th diagonal entry is $m_i(t)$, and
   $F(t) =D^{-1}_{\mathbb{N}(t)}A'_{\mathbb{N}(t)}$.  Note that $F(t)$
  is  a stochastic matrix and sometimes referred to as a {\em flocking matrix} in the literature \cite{vicsekmodel}.
It is straightforward to verify that
\eq{y(t+1) =P(F(t)\otimes I)Py(t),\;\;\;\;\; t\geq 1\label{sys}}
where $\otimes$ denotes the Kronecker product, $I$ denotes the $n\times n$ identity matrix,
and
$P = {\rm diagonal}$ $\{P_1,P_2,\ldots, P_m\}$ is an $mn\times mn$ block diagonal matrix.
  Note that $P^2=P$ because each $P_i$ is idempotent.  We will use this fact without special mention in the sequel.
It is also worth noting that $\bbb{N}(t)=\gamma(F(t))$ where
$\gamma(F(t))$ is the graph of $F(t)$ whose definition has been given in \S \ref{notion}.

\section{$\scr{D}$-Connectivity} \label{pconnect}

We begin with the special case in which the linear equation $Ax=b$ has a unique solution
and the neighbor graph $\bbb{N}(t)$ is fixed and independent of time.
Our reason for considering this special case first will become clear shortly.

Note that $Ax=b$ has a unique solution exactly when $\bigcap_{i=1}^m \ker A_i = \{0\}$.
Since $\ker A_i  = \scr{P}_i$, $i\in\mathbf{m}$,
where $\scr{P}_i$ denotes the column span of $P_i$,
the uniqueness assumption is equivalent to the condition
\eq{\bigcap_{i=1}^m\scr{P}_i = \{0\}\label{assmp}}
In this case, all $x_i(t)$ approach the unique solution in the limit as $t\rightarrow\infty$
if and only if $y(t)\rightarrow 0$.
If we further assume that $\bbb{N}(t)=\bbb{N}$ is independent of time,
the linear recursion equation \rep{sys} simplifies to a linear time-invariant system
$y(t+1) =P(F\otimes I)Py(t)$ where $F$ is the flocking matrix of $\bbb{N}$.
Then, the algorithm \rep{a1} causes all $x_i$ to converge to the unique solution of $Ax=b$ exponentially fast
if and only if $P(F\otimes I)P$ is a discrete-time
stability matrix.\footnote{
A square matrix $M$ is called a {\em discrete-time stability matrix} if
the largest magnitude of all eigenvalues of $M$ is strictly less than $1$.
}
We are particularly interested in the weakest possible graph connectivity condition for $P(F\otimes I)P$ to be a discrete-time
stability matrix.
To state the condition, we need the following concepts.

We write $\scr{V} \triangleq \mathbf{m}$ to denote the vertex set.
Suppose that $\scr{E}$ is a nonempty subset of $\scr{V}$.
It follows immediately that $\bigcap_{i\in\scr{E}} \scr P_i \supset \bigcap_{i\in\scr{V}} \scr P_i$.
We say that $\scr{E}$ is a {\em fully populated} set  if
$\bigcap_{i\in\scr{E}} \scr P_i = \bigcap_{i\in\scr{V}} \scr P_i$.
We say that $\scr{E}$ is a {\em partially populated} set
if $\bigcap_{i\in\scr{V}} \scr P_i$ is a proper subset of $\bigcap_{i\in\scr{E}} \scr P_i$.
Clearly a partially populated set $\scr{E}$ must be a nonempty proper subset of $\scr{V}$.

From the above definitions, it can be seen that a subset of $k$ agents can solve the original linear equation $Ax=b$
without any information from the remaining agents
if and only if the set of their labels is fully populated.
Note that in the case when $Ax=b$ has a unique solution,
$\scr{E}$ is fully populated if $\bigcap_{i\in\scr{E}} \scr P_i=\{0\}$ and
is partially populated if $\bigcap_{i\in\scr{E}} \scr P_i$ is a nonzero subspace.

Let $\bbb{G}$ be a directed graph with vertex set $\scr{V}$.
It is natural to call a vertex $j$ a {\em neighbor} of vertex $i$
if $(j,i)$ is an arc in $\bbb{G}$.
We say that vertex $i\in\scr{V}$ has a {\em neighbor} in
$\scr{W}\subset\scr{V}$ if there is a vertex $j\in\scr{W}$ which is a neighbor of $i$.
We say that a set $\scr U\subset\scr V$ has a {\em neighbor} in
$\scr{W}\subset\scr{V}$ if there exist vertices $i\in\scr{U}$ and
$j\in\scr{W}$ such that $j$ is a neighbor of $i$.

\begin{proposition}
Suppose that \rep{assmp} holds.
Then, $P(F\otimes I)P$ is a discrete-time stability matrix
if and only if every partially populated
subset $\scr{E}\subset \scr{V}$ in $\gamma(F)$ has at least one neighbor in $\scr{V}\setminus\scr{E}$,
where $\scr{V}\setminus\scr{E}$ denotes the complement of $\scr{E}$ in $\scr{V}$.
\label{ns}\end{proposition}
The proof of Proposition \ref{ns} can be found in the appendix.


The necessary and sufficient graph-theoretic condition in Proposition \ref{ns} can be interpreted as follows:
{\em Any subset of agents need information from outside the subset if they cannot solve the problem by themselves}.
Prompted by this, we define a new notion of graph connectivity as follows.

Let $\bbb{G}$ be a directed graph with vertex set $\scr{V}$
and let $2^{\scr{V}}$ denote the power set of $\scr{V}$.
We say that a collection of nonempty proper subsets $\scr{W}\subset 2^{\scr{V}}$ is {\em connected} by $\bbb{G}$ if
each subset $\scr{E}\in\scr{W}$ has at least one neighbor in $\scr{V}\setminus\scr{E}$.

Let $\scr{D}\subset 2^{\scr{V}}$ be the collection of all partially populated subsets of $\scr{V}$.
We say that a directed graph $\bbb{G}$ is {\em $\scr{D}$-connected} if
$\scr{D}$ is connected by $\bbb{G}$.

Note that the notion of $\scr{D}$-connectivity depends on $\scr{D}$
and thus on the set $\{P_1,P_2,\ldots,P_m\}$.
{\color{black}
Therefore, the notion is dependent on the data about $\matt{A & b}$
available to the agents.
In some applications \cite{estimation07,calibration}, the data about $\matt{A & b}$
may be affected by noise, in which case the notion of $\scr{D}$-connectivity
also depends on the noise. In this paper, we assume that the data about $\matt{A & b}$
is not corrupted by any noise, i.e., each agent $i$ accurately knows $A_i$ and $b_i$;
thus, the notion of $\scr{D}$-connectivity is completely determined by the data about $\matt{A & b}$
available to the agents.
}
While it is easy to see that every strongly connected graph is $\scr{D}$-connected no matter what $\scr{D}$ is,
the following example shows that with respect to some $\scr{D}$,
there are $\scr{D}$-connected graphs which are not strongly connected.

\begin{figure}[h]
\centering
\includegraphics[height=70mm]{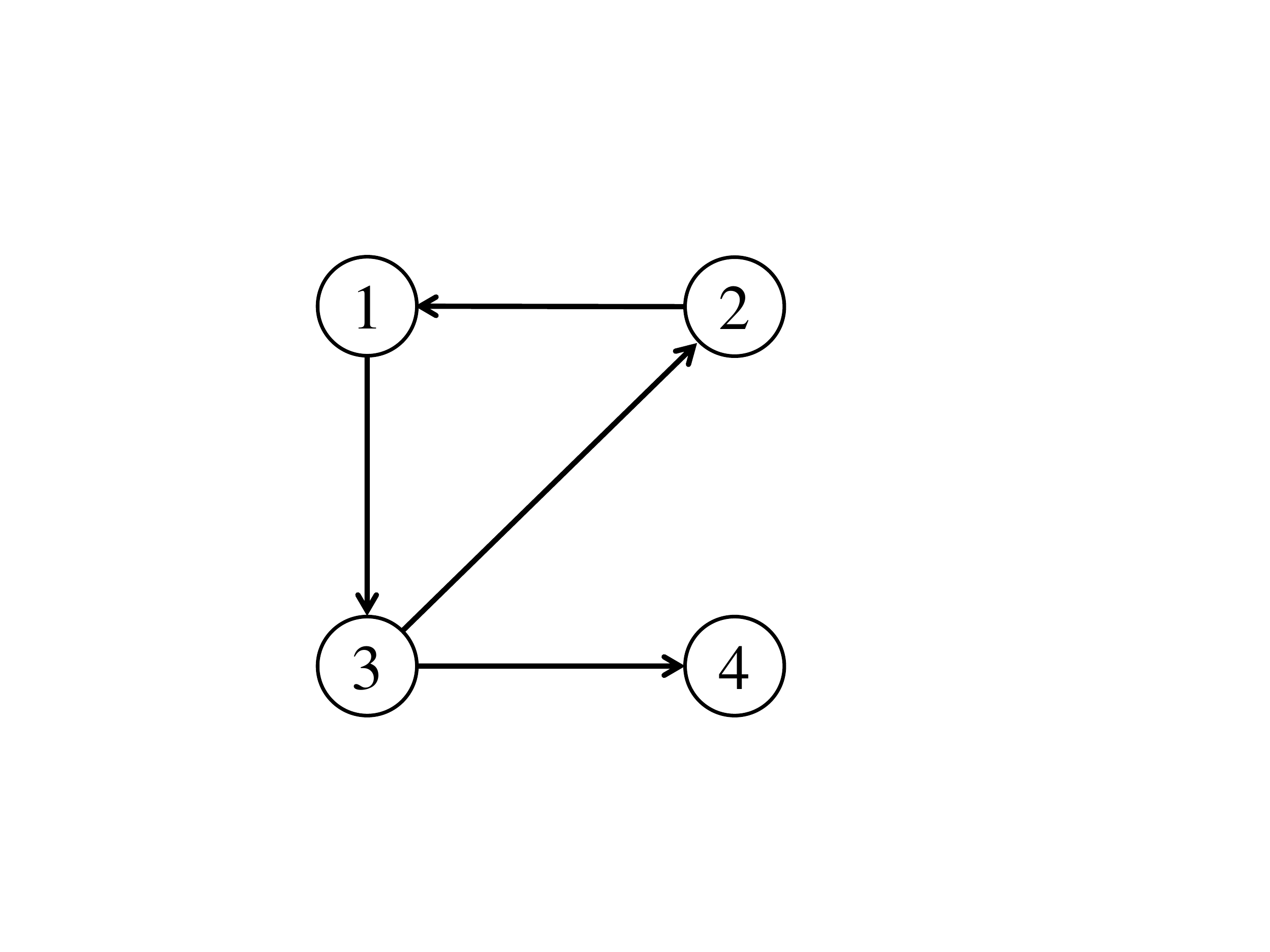}
\caption{A $\scr{D}$-connected graph which is not strongly connected}
\label{1}
\end{figure}

Consider a network consisting of $4$ agents. Suppose that
$$A=\matt{1 & 2 & 3 &4 \cr 5 & 6& 7& 8 \cr 3& 6 & 9 & 12 \cr 10 & 12 & 14 &16}$$
and that each agent $i$ knows the $i$th row of $A$, $i\in\{1,2,3,4\}$.
Since $A_3=3A_1$ and $A_4=2A_2$, there hold $P_1=P_3$ and $P_2=P_4$.
Thus, in this case, all partially populated subsets are
$\{1\}$, $\{2\}$, $\{3\}$, $\{4\}$, $\{1,3\}$ and $\{2,4\}$, i.e.,
$\scr{D}=\{\{1\}, \{2\}, \{3\}, \{4\}, \{1,3\}, \{2,4\}\}$.
It is straightforward to verify that the graph in Figure \ref{1} is $\scr{D}$-connected, but not strongly connected.

Thus, the set of strongly connected graphs is a subset of the set of $\scr{D}$-connected graphs no matter what $\scr{D}$ is,
and can be a proper subset of the set of $\scr{D}$-connected graphs, depending on $\scr{D}$.
A sufficient condition under which the two properties of being strongly connected and
$\scr{D}$-connected are equivalent, is discussed in \S \ref{equal}.

\begin{figure}[h]
\centering
\includegraphics[height=70mm]{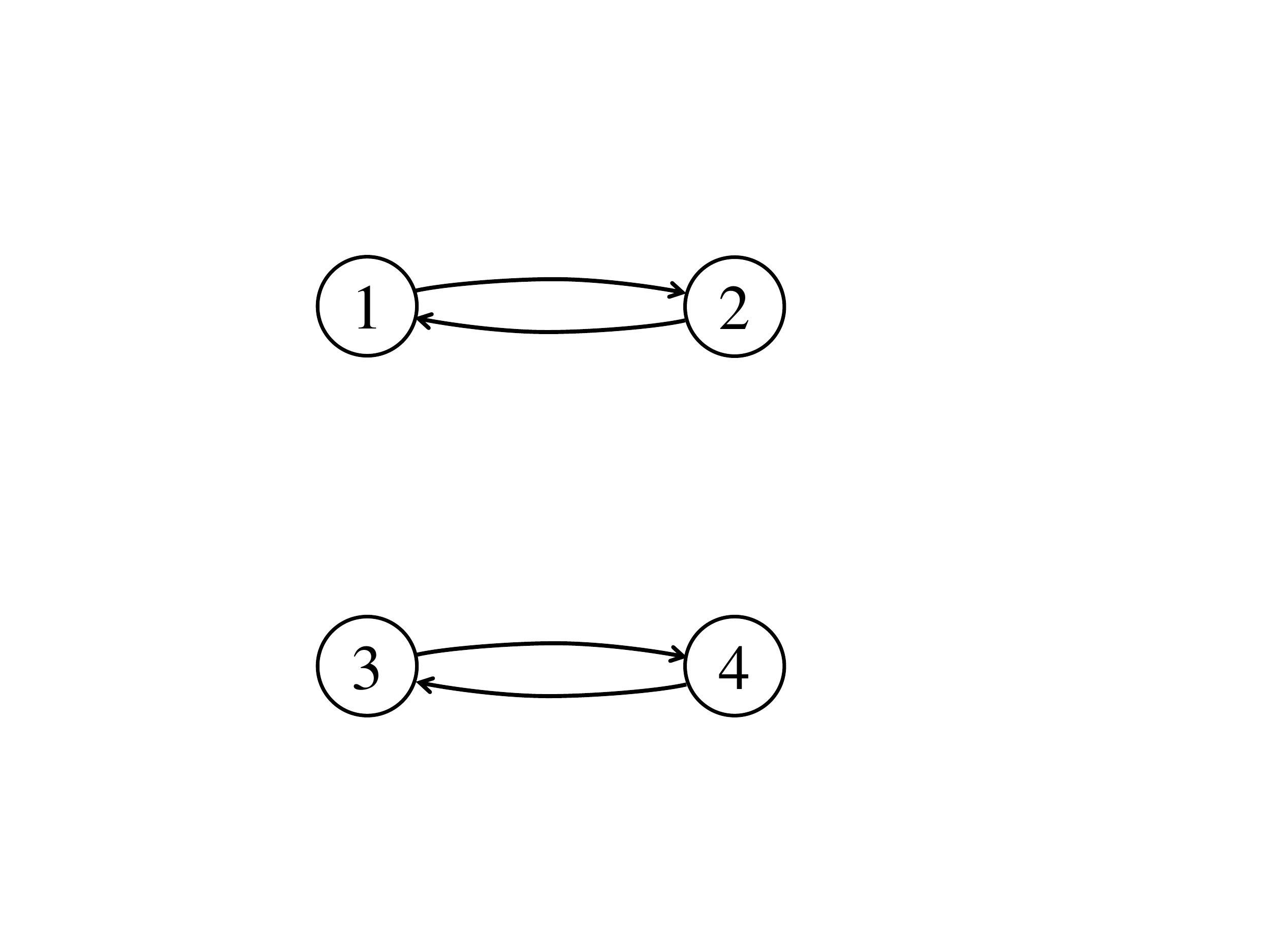}
\caption{A $\scr{D}$-connected graph which is not rooted}
\label{3}
\end{figure}

\begin{figure}[h]
\centering
\includegraphics[height=70mm]{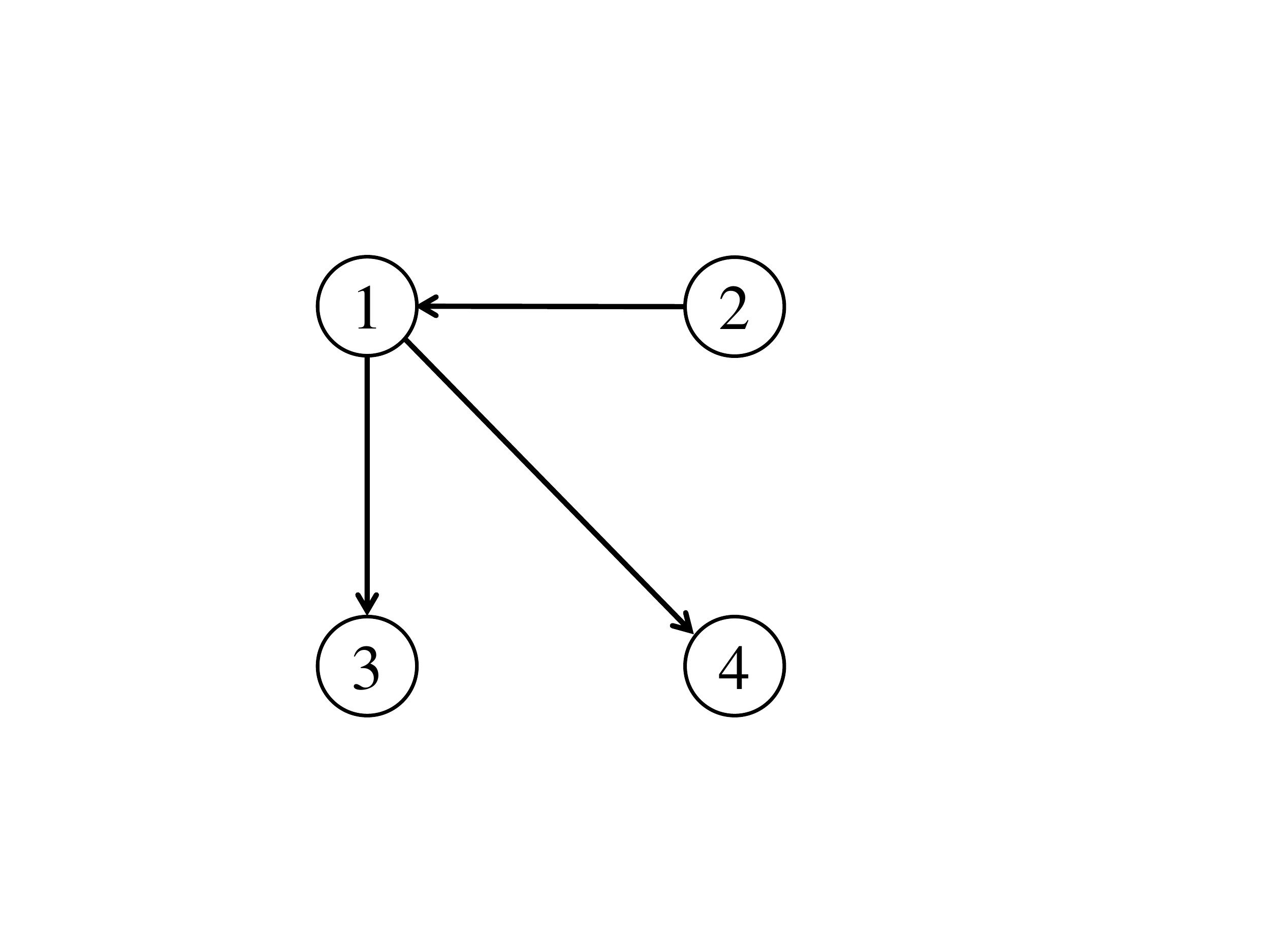}
\caption{A rooted graph which is not $\scr{D}$-connected}
\label{2}
\end{figure}

It is worth noting that the notions of $\scr{D}$-connectivity and root connectivity are not comparable.
In particular, there are $\scr{D}$-connected graphs which are not rooted,
and vice versa.
For instance, in the preceding $4$-agent network,
it is straightforward to verify that the graph in Figure \ref{3} is $\scr{D}$-connected but not rooted,
and that the graph in Figure \ref{2} is rooted but not $\scr{D}$-connected.

{\color{black}
It is also worth emphasizing that whether the graphs in Figures \ref{1}-\ref{2}
are $\scr{D}$-connected or not, depends on the matrix $A$ and how its rows are partitioned.
}

Proposition \ref{ns} implies that the notion of $\scr{D}$-connectivity, which is less restrictive
than strong connectivity, is the weakest possible graph connectivity condition for
 exponential convergence of the algorithm \rep{a1} in the special case when
$\bbb{N}(t)$ is fixed and $Ax=b$ has a unique solution.
In the next section, we will show that this connectivity notion is also
appropriate to the analysis of general cases.

\section{Main Results} \label{main}

We now turn to the general cases in which the neighbor graph $\bbb{N}(t)$ may change over time.
We consider the cases of unique solution and nonunique solution separately.

To state our main results, we need some naturally extended notions of graph connectivity,
as was done for strong and root connectivity in \S \ref{notion}.
We say that a finite sequence of directed graphs $\bbb{G}_1,\bbb{G}_2,\ldots,\bbb{G}_p$ with the same vertex set
is {\em jointly $\scr{D}$-connected} if the composition
$\mathbb{G}_{q}\circ\mathbb{G}_{q-1}\circ  \cdots\circ\mathbb{G}_{1}$ is $\scr{D}$-connected.
We say that an infinite  sequence of directed graphs $\mathbb{G}_1, \mathbb{G}_2,\ldots$ with the same vertex set
is {\em repeatedly jointly $\scr{D}$-connected}
if there exist finite positive integers $l$ and $\tau_0$ such that for any integer $k\geq 0$,
the finite sequence  $\mathbb{G}_{\tau_0+kl},\mathbb{G}_{\tau_0+kl+1},\ldots,\mathbb{G}_{\tau_0+(k+1)l-1}$
is jointly $\scr{D}$-connected.
If such integers $l$ and $\tau_0$ exist, we say that $\mathbb{G}_{\tau_0}, \mathbb{G}_{\tau_0+1},\ldots$
is repeatedly jointly $\scr{D}$-connected {\em by subsequences of length $l$}.

The following theorem gives a necessary and sufficient graph-theoretic condition
for the algorithm \rep{a1} to converge exponentially fast in the case when $Ax=b$ has a unique solution.

\begin{theorem}
Suppose that $Ax=b$ has a unique solution and that each agent $i$ updates its state $x_i(t)$
according to algorithm \rep{a1}. Then, there exists a nonnegative constant $\lambda<1$
for which all $x_i(t)$ converge to the unique solution to $Ax=b$ as $t\rightarrow\infty$
at the same rate as $\lambda^t$ converges to $0$ if and only if
the sequence of neighbor graphs $\bbb{N}(1),\bbb{N}(2),\ldots$ is repeatedly jointly $\scr{D}$-connected.
\label{main1}\end{theorem}

The proof of Theorem \ref{main1} can be found in \S \ref{unique}.
{\color{black}
It is worth emphasizing that compared with Theorem 2 in \cite{tac},
the above theorem relaxes strong connectivity to $\scr D$-connectivity without the non-redundancy assumption.
}

Using the same arguments as in the proof of Corollary 1 in \cite{tac},
we have the following result on the convergence rate.

{\color{black}
\begin{corollary}
Suppose that $Ax = b$ has a unique solution $x^*$.
Suppose that each neighbor graph $\bbb N(t)$, $t \ge 1$, is $\scr D$-connected;
let $\scr F$ be the set of all flocking matrices corresponding to $\bbb N(t)$, $t \ge 1$.
Let
$$\lambda = \left(1-\frac{(m-1)(1-\rho)}{m^\tau}\right)^{\frac{1}{\tau}}$$
where $\tau$
is a finite positive integer, not depending on $\scr F$, such that
for any  $p\geq \tau$,  the matrix
$P(F_{p}\otimes I)P(F_{p-1}\otimes I)\cdots  P(F_1\otimes I)P$, with each $F_i\in\scr F$, $i\in\{1,2,\ldots,p\}$,
  is a contraction in the mixed  matrix norm,\footnote{
Such a $\tau$ exists as shown in Proposition \ref{p}.
}
and
$$\rho = \max_{\scr C} |P_{j_1}P_{j_2}\cdots P_{j_{\tau+1}}|_2$$
in which $|\cdot|_2$ denotes the induced two-norm, and
$\scr C$ is the set of all products of the matrices in $\{P_1,P_2,\ldots,P_m\}$ of length $\tau+1$
such that each matrix in $\{P_1,P_2,\ldots,P_m\}$ occurs in the product at least once.\footnote{
The set $\scr{C}$ is compact and $\rho$ is less than $1$.}
Then, all $x_i(t)$ converge to $x^*$ as $t \rightarrow 0$ as fast as $\lambda^t$ converges
to $0$.
\end{corollary}
}

This result can readily be extended to the case when $Ax = b$ has more
than one solution.

For the case when $Ax=b$ has more than one solution, the necessary and sufficient graph-theoretic condition
has a different characterization than it does for the uniqueness case.
In the special case when $A=0$ (and consequently $b=0$), the problem reduces to
an unconstrained consensus problem in which
all $m$ states converge exponentially fast to the same value if and only if
the sequence of neighbor graphs $\bbb{N}(1),\bbb{N}(2),\ldots$ is repeatedly jointly rooted
\cite{luc,cdc14}.
The following theorem gives a necessary and sufficient graph-theoretic condition
for the algorithm \rep{a1} to converge exponentially fast in the case when $A\neq 0$.

\begin{theorem}
Suppose that $Ax=b$ has more than one solution, $A\neq 0$, and that each agent $i$ updates its state $x_i(t)$
according to algorithm \rep{a1}. Then, there exists a nonnegative constant $\lambda<1$
for which all $x_i(t)$ converge to the same solution to $Ax=b$ as $t\rightarrow\infty$
at the same rate as $\lambda^t$ converges to $0$ if and only if
the sequence of neighbor graphs $\bbb{N}(1),\bbb{N}(2),\ldots$ is
repeatedly jointly rooted and $\scr{D}$-connected.
\label{main2}\end{theorem}
The proof of Theorem \ref{main2} can be found in \S \ref{nonunique}.

As noted earlier, the set of strongly connected graphs can be a proper subset of
the set of $\scr{D}$-connected graphs and is a proper subset of the set of rooted graphs.
Thus, the graph-theoretic conditions given in Theorem \ref{main1} and Theorem \ref{main2}
are both less restrictive than the repeatedly jointly strongly connected condition
required in \cite{cdc13.1}.
As the next section will demonstrate, the proofs of these theorems are more complicated
and challenging than those in \cite{cdc13.1}.

\section{Analysis} \label{analysis}

The aim of this section is to give proofs of Theorem \ref{main1} and Theorem \ref{main2}.
We begin with the case when $Ax=b$ has a unique solution (Theorem \ref{main1}).

\subsection{Uniqueness} \label{unique}

We first prove the necessity of Theorem \ref{main1}.
For this, we need the following lemma.

\begin{lemma}
Suppose that \rep{assmp} holds. Let $S_1,S_2,\ldots, S_p$  be a finite sequence of $m\times m$ stochastic matrices
with positive diagonal entries.
If $\gamma(S_1),\gamma(S_2),\ldots,\gamma(S_p)$ is not jointly $\scr{D}$-connected, then the matrix
$P(S_{p}\otimes I)P(S_{p-1}\otimes I)\cdots  P(S_1\otimes I)P$
has an eigenvalue at $1$.
\label{eig1}\end{lemma}

\noindent
{\bf Proof of Lemma \ref{eig1}:}
Since $\gamma(S_1),\gamma(S_2),\ldots,\gamma(S_p)$ is not jointly $\scr{D}$-connected, the composed graph
$$\mathbb{G} = \gamma(S_p)\circ \gamma(S_{p-1})\circ \cdots\circ  \gamma(S_1)$$
is not $\scr{D}$-connected.
Then, there exists a partially populated subset $\scr{E}$, which is a proper subset of $\scr{V}=\mathbf{m}$,
such that $\scr{E}$ does not have a neighbor in $\scr{V}\setminus\scr{E}$.
Let $\scr{E} =\{i_1,i_2,\ldots,i_k\}$, $k<m$.
Since each $\gamma(S_i)$, $i\in\mathbf{p}$, has self-arcs at all vertices,
the arcs of each $\gamma(S_i)$ are arcs of $\bbb{G}$.
Since there is no arc from $\scr{V}\setminus\scr{E}$ to $\scr{E}$ in $\bbb{G}$,
there is no arc from $\scr{V}\setminus\scr{E}$ to $\scr{E}$ in each $\gamma(S_i)$, $i\in\mathbf{p}$.
Let $\pi$ be any permutation on $\scr{V}$ for which
$\pi(i_j) = j$, $j\in\mathbf{k}$, and let $Q$ be the corresponding permutation matrix.
Then, for each $i\in\mathbf{p}$, the transformation $S_i\rightarrow QS_iQ'$ block triangularizes $S_i$.
Set $\bar Q = Q\otimes I$. Note that $\bar Q$ is a permutation matrix and
that $\bar QP\bar Q'$ is a block diagonal, orthogonal projection matrix
where the $j$th diagonal block is $P_{\pi(i_j)}$, $j\in\mathbf{k}$.
Since each $QS_iQ'$ is block triangular, so are the matrices $\bar QP(S_i\otimes I)P\bar Q'$, $i\in\mathbf{p}$.
Thus, for each $i\in\mathbf{p}$, there are matrices $A_i$, $B_i$ and $C_i$ such that
 $$\bar QP(S_i\otimes I)P\bar Q' = \matt{A_i & 0\cr B_i &C_i}$$

For each $i\in\mathbf{p}$, let $\bar S_i$ be that $k\times k$ submatrix of $S_i$ whose $pq$th entry is
the $i_pi_q$th entry of $S_i$ for all $p,q\in\mathbf{k}$.
In other words, $\bar S_i$ is that submatrix of $S_i$ obtained by deleting rows and
columns whose indices are not in $\scr{E}$.
Since $S_i$ is a stochastic matrix and there are no arcs from $\scr{V}\setminus\scr{E}$ to $\scr{E}$,
it follows that $\bar S_i$ is a stochastic matrix and
its graph $\gamma(\bar S_i)$ is the subgraph of $\gamma(S_i)$ induced by $\scr{E}$.
Set $P_{\scr{E}} = {\rm diagonal}\{P_{i_1},P_{i_2},\ldots,P_{i_k}\}$.
Then, it is straightforward to verify that
$$A_i = P_{\scr{E}} \displaystyle\left(\bar S_i\otimes I\right)P_{\scr{E}}$$
Since $\scr{E}$  is a partially populated subset, it follows by definition that
$\bigcap_{j\in\scr{E}} \scr P_j \neq \{0\}$. Then, for any nonzero vector $z$ in $\bigcap_{j\in\scr{E}} \scr P_j $, there holds
$A_i \bar z = \bar z$,
where $\bar z= {\rm stack} \{z,z,\ldots,z\}$.
 Note that
\begin{eqnarray*}
&&\bar Q\displaystyle\left(P(S_{p}\otimes I)P(S_{p-1}\otimes I)\cdots  P(S_1\otimes I)P\right)\bar Q' \\&=&
\displaystyle\left(\bar QP(S_{p}\otimes I)P\bar Q'\right) \cdots  \displaystyle\left(\bar QP(S_1\otimes I)P\bar Q'\right)\\
&=& \matt{A & 0\cr B &C}
\end{eqnarray*}
where $A=A_pA_{p-1}\cdots A_1$. It follows that $A\bar{z} = \bar{z}$,
 so $A$ has an eigenvalue at $1$. Therefore, the matrix
 $P(S_{p}\otimes I)P(S_{p-1}\otimes I)\cdots  P(S_1\otimes I)P$
has an eigenvalue at $1$.
\hfill$\qed$

To proceed, we need a special ``mixed matrix norm'' introduced in \cite{cdc13.1}.
Let    $|\cdot |_{\infty}$ denote   the
  induced infinity  norm   and write
 $\R^{mn\times mn}$ for  the vector space of all $m\times m$  block matrices $Q = \matt{Q_{ij}}$
whose $ij$th entry     is an  $n\times n$ matrix $Q_{ij}\in\R^{n\times n}$.
  We define the {\em mixed matrix norm} of $Q\in\R^{mn\times mn}$, written $||Q||$, to be
$$||Q|| = |\langle Q\rangle |_{\infty}$$
 where $\langle Q\rangle $ is the $m\times m$ matrix in $\R^{m\times m}$  whose $ij$th entry is
{\color{black} $|Q_{ij}|_2$,
 where $|\cdot|_2$ denotes the induced two-norm.}
It has been shown in \cite{cdc13.1} that $||\cdot ||$ is a sub-multiplicative norm (see Lemma 3 of \cite{cdc13.1}).

{\color{black}
For the matrices of the form $P(F(t)\otimes I)P$ defined in \rep{sys}, more can be said.
It has been shown in \cite{tac} that such matrices are non-expansive in the mixed matrix norm,
i.e., $\|P(F(t)\otimes I)P\|\le 1$
(see Proposition 1 in \cite{tac}).
Since $||\cdot ||$ is a sub-multiplicative norm,
there holds $\|\Phi(t,\tau)\|\le 1$ for all $t$ and $\tau$,
where $\Phi(t,\tau)$ denotes the state transition matrix of $P(F(t)\otimes I)P$.
}

\noindent
{\bf Proof of Theorem \ref{main1} (Necessity):}
In the case when $Ax=b$ has a unique solution, all $x_i(t)$ in \rep{a1}
converge to the unique solution exponentially fast precisely when
the linear system \rep{sys} is exponentially stable.
Since exponential stability and uniform asymptotic stability are equivalent properties for linear systems,
it will be sufficient to show that uniform asymptotic stability of \rep{sys} implies that
the sequence of neighbor graphs $\bbb{N}(1),\bbb{N}(2),\ldots$ is repeatedly jointly $\scr{D}$-connected.
Suppose therefore that the system \rep{sys} is uniformly asymptotically stable.

To establish the claim, suppose that, to the contrary, $\bbb{N}(1),\bbb{N}(2),\ldots$
is not repeatedly jointly $\scr{D}$-connected. Then,
the negation of the definition of repeatedly jointly $\scr{D}$-connected graphs implies that
for any pair of positive integers $p$ and $q$, there is an integer $k\ge q$ such that
the composed graph
$\mathbb{N}(k+p-1)\circ  \cdots\circ \mathbb{N}(k+1)\circ \mathbb{N}(k)$
is not $\scr{D}$-connected.

Let $\Phi(t,\tau)$ be the state transition matrix of $P(F(t)\otimes I)P$.
Since \rep{sys} is uniformly asymptotically stable, for each real number $e>0$, there exist
positive integers $t_e$ and $T_e$ such that $||\Phi(t+T_e,t)||<e$ for all $t\ge t_e$.
Set $e=1$ and let $t_1$ and $T_1$ be any pair of such integers.
It follows from the preceding that there is an integer $t_2 \ge t_1$ such that
the composed graph
$$\mathbb{N}(t_2+T_1-1)\circ  \cdots\circ \mathbb{N}(t_2+1)\circ \mathbb{N}(t_2)$$
is not $\scr{D}$-connected.
Since $t_2 \ge t_1$, the hypothesis of uniform asymptotic stability ensures that
\eq{\|\Phi(t_2+T_1,t_2)\|<1\label{x}}
Note that $\bbb{N}(t)=\gamma(F(t))$ for all $t\in\{t_2,t_2+1,\ldots,t_2+T_1-1\}$ and
$$\Phi(t_2+T_1,t_2) = P(F(t_2+T_1-1)\otimes I)P\cdots (F(t_2)\otimes I)P$$
By Lemma \ref{eig1},  $\Phi(t_2+T_1,t_2)$
has an eigenvalue at $1$, so $\|\Phi(t_2+T_1,t_2)\|=1$.
But this contradicts \rep{x}. Therefore, $\bbb{N}(1),\bbb{N}(2),\ldots$
is repeatedly jointly $\scr{D}$-connected.
\hfill$\qed$

We now turn to the proof of sufficiency.
The sufficiency of Theorem \ref{main1} is a consequence of the  following result.

\begin{proposition}
Suppose that \rep{assmp} holds. Let $S_1,S_2,\ldots$  be a sequence of $m\times m$ stochastic matrices
whose corresponding sequence of graphs $\gamma(S_1),\gamma(S_2),\ldots$ is repeatedly jointly $\scr{D}$-connected
by subsequences of length $l$.
Then, there is a finite positive integer $\tau$, not depending on $P$, such that
for any  $p\geq \tau$,  the matrix
$P(S_{p}\otimes I)P(S_{p-1}\otimes I)\cdots  P(S_1\otimes I)P$
  is a contraction in the mixed  matrix norm.
\label{p}\end{proposition}

To prove Proposition \ref{p}, we need a few concepts.

We call a vertex $i$ in a directed graph $\mathbb{G}$  a {\em sink}  of $\mathbb{G}$
  if for each other vertex $j$  of $\mathbb{G}$, there is a directed  path from $j$ to $i$.
 We say that
$\mathbb{G}$ is {\em sunk  at $i$}  if $i$ is in fact  a sink. Thus,
$\mathbb{G}$ is sunk  at $i$ whenever $i$  is {\em reachable} from each other vertex of $\mathbb{G}$
 along a directed path within the graph.
 $\mathbb{G}$ is  {\em strongly sunk at $i$} if $i$ is  reachable from  each other vertex of $\mathbb{G}$
  along a  directed path of length $1$. Thus,
$\mathbb{G}$ is strongly sunk at $i$ if $i$ is an observer\footnote{
We refer to $k$ as an {\em observer} of $i$ if $(i,k)$ is an arc in $\bbb{G}$.}
of every other vertex in the graph.
 By a {\em sunk graph} $\mathbb{G}$ is meant a directed graph which possesses at least one sink.
A {\em strongly sunk graph} is a graph which has at least  one vertex at which it is strongly sunk.
It is worth noting that a directed graph $\bbb{G}$ is (strongly) sunk if
its dual graph $\bbb{G}'$ is (strongly) rooted.\footnote{
The {\em dual graph} of a directed graph $\bbb{G}$ is that
graph which results when the arcs in $\bbb{G}$ are reversed.}

To proceed, we need a concept from \cite{reachingp1}.
Let $\scr{G}$ denote the set of all directed graphs with $m$ vertices.
By the {\em neighbor function} of a directed graph $\bbb{G}\in\scr{G}$ with the vertex set $\scr{V}=\mathbf{m}$,
denoted by $\beta(\bbb{G},\cdot)$,
we mean the function $\beta(\bbb{G},\cdot): 2^{\scr{V}}\rightarrow 2^{\scr{V}}$
which assigns to each subset $\scr{S}\subset\scr{V}$, the subset of vertices in $\scr{V}$
which are neighbors of $\scr{S}$ in $\bbb{G}$.
Thus, $j\in\beta(\bbb{G},i)$ whenever $(j,i)$ is an arc in $\bbb{G}$.
Note that if $\bbb{G}_q\in\scr{G}$ and $\bbb{G}_p\in\scr{G}_{sa}$, then
\eq{\beta(\bbb{G}_q,\scr{S})\subset \beta(\bbb{G}_q\circ\bbb{G}_p,\scr{S}), \hspace{.3in}\scr{S}\in 2^{\scr{V}}\label{2chain}}
since $\bbb{G}_p\in\scr{G}_{sa}$ implies that the arcs in $\bbb{G}_q$ are all arcs in $\bbb{G}_q\circ\bbb{G}_p$.
Neighbor functions have the following important property.

\begin{lemma}
For all $\bbb{G}_p,\bbb{G}_q\in\scr{G}$ and any nonempty subset $\scr{S}\subset\scr{V}$, there holds
$$\beta\displaystyle\left(\bbb{G}_p, \beta(\bbb{G}_q,\scr{S})\right) = \beta\displaystyle\left(\bbb{G}_q\circ\bbb{G}_p,\scr{S}\right)$$
\label{func1}\end{lemma}

\noindent
{\bf Proof of Lemma \ref{func1}:}
We first show that $\beta(\bbb{G}_p, \beta(\bbb{G}_q,\scr{S}))$ $\subset \beta(\bbb{G}_q\circ\bbb{G}_p,\scr{S})$.
Suppose that $i\in\beta(\bbb{G}_p, \beta(\bbb{G}_q,\scr{S}))$. Then, $(i,j)$ is an arc in $\bbb{G}_p$
for some $j\in\beta(\bbb{G}_q,\scr{S})$.
Hence, $(j,k)$ is an arc in $\bbb{G}_q$ for some $k\in\scr{S}$. In view of the definition of composition,
$(i,k)$ is an arc in $\bbb{G}_q\circ \bbb{G}_p$, so $i\in\beta(\bbb{G}_q \circ \bbb{G}_p,\scr{S})$.
Since this holds for all $i\in\scr{V}$, it follows that
$\beta(\bbb{G}_p, \beta(\bbb{G}_q,\scr{S})) \subset \beta(\bbb{G}_q\circ\bbb{G}_p,\scr{S})$.

For the reverse inclusion, suppose that $i\in\beta(\bbb{G}_q\circ\bbb{G}_p,\scr{S})$,
which implies that $(i,k)$ is an arc in $\bbb{G}_q \circ \bbb{G}_p$ for some $k\in\scr{S}$.
By the definition of composition, there exists an $j\in\scr{V}$ such that $(i,j)$ is an arc in $\bbb{G}_p$
and $(j,k)$ is an arc in $\bbb{G}_q$. Then, $j\in\beta(\bbb{G}_q,\scr{S})$.
Thus, $i\in\beta(\bbb{G}_p,\beta(\bbb{G}_q,\scr{S}))$. Since this holds for all $i\in\scr{V}$,
it follows that
$\beta(\bbb{G}_q\circ\bbb{G}_p,\scr{S})\subset\beta(\bbb{G}_p, \beta(\bbb{G}_q,\scr{S}))$.
Therefore, the lemma is true.
\hfill$\qed$

Let us note that each subset $\scr{S}\subset\scr{V}$ induces a unique subgraph of $\bbb{G}$
with vertex set $\scr{S}$ and arc set $\scr{A}$ consisting of those arcs $(i,j)$ of $\bbb{G}$ for which
both $i$ and $j$ are vertices of $\scr{S}$. This together with the natural partial ordering of $\scr{V}$
by inclusion provides a corresponding partial ordering of $\bbb{G}$.
Thus, if $\scr{S}_1$ and $\scr{S}_2$ are subsets of $\scr{V}$ and $\scr{S}_1\subset\scr{S}_2$,
then $\bbb{G}_1\subset\bbb{G}_2$ where for $i\in\{1,2\}$, $\bbb{G}_i$ is the subgraph of $\bbb{G}$ induced
by $\scr{S}_i$. For any $v\in\scr{V}$, there is a unique largest subgraph sunk at $v$, namely
the graph induced by the vertex set
$\scr{V}(v) = \{v\}\cup\beta(\bbb{G},v)\cup \cdots \cup \beta^{m-1}(\bbb{G},v)$
where $\beta^i(\bbb{G},\cdot)$ denotes the composition of $\beta(\bbb{G},\cdot)$ with itself $i$ times.
We call this graph, the {\em sunk graph generated by $v$}.
Note that $\scr{V}(v)$ is the smallest $\beta(\bbb{G},\cdot)$ - invariant subset of $\scr{V}$
which contains $v$.
The sunk graph generated by a vertex of a $\scr{D}$-connected graph
has the following property.

\begin{lemma}
Suppose that $\bbb{G}$ is a graph  in $\scr{G}_{sa}$ which is $\scr{D}$-connected.
Then, for each vertex $v$ of $\bbb{G}$, there holds
$$\bigcap_{i\in\scr{V}(v)} \scr P_i = \{0\}$$
where
$\scr{V}(v) = \{v\}\cup\beta(\bbb{G},v)\cup \cdots \cup \beta^{m-1}(\bbb{G},v)$.
\label{easy}\end{lemma}

\noindent
{\bf Proof of Lemma \ref{easy}:}
To prove the lemma, suppose that, to the contrary,
$\bigcap_{i\in\scr{V}(v)} \scr P_i$ is a nonzero subspace.
Then, $\scr{V}(v)$ is a partially populated subset of $\scr{V}$.
Since $\scr{V}(v)$ is the vertex set of the sunk graph generated by $v$,
it follows that $\scr{V}(v)$ has no neighbor in $\scr{V}\setminus\scr{V}(v)$.
But this is impossible since $\bbb{G}$ is $\scr{D}$-connected.
Thus, the lemma is true.
\hfill$\qed$

The proof of Proposition \ref{p} depends on the following lemmas.

\begin{lemma}
Let $\bbb{G}_p$ and $\bbb{G}_q$ be graphs in $\scr{G}_{sa}$. If $\bbb{G}_p$ is sunk at $v$
and $\beta(\bbb{G}_q,v)$ is a strictly proper subset of $\scr{V}$,
then $\beta(\bbb{G}_q,v)$ is also a strictly proper subset of $\beta(\bbb{G}_q \circ \bbb{G}_p,v)$.
\label{func2}\end{lemma}

\noindent
{\bf Proof of Lemma \ref{func2}:}
Note that $\beta(\bbb{G}_q,v)\subset \beta(\bbb{G}_q\circ\bbb{G}_p,v)$
because of \rep{2chain}.
Thus, if $\beta(\bbb{G}_q,v)$ is not a strictly proper subset of $\beta(\bbb{G}_q\circ\bbb{G}_p,v)$,
then $\beta(\bbb{G}_q,v)=\beta(\bbb{G}_q\circ\bbb{G}_p,v)$, so
$\beta(\bbb{G}_q\circ\bbb{G}_p,v)\subset\beta(\bbb{G}_q,v)$.
In view of Lemma \ref{func1},
$\beta(\bbb{G}_q\circ\bbb{G}_p,v)=\beta(\bbb{G}_p, \beta(\bbb{G}_q,v))$.
Thus, $\beta(\bbb{G}_p, \beta(\bbb{G}_q,v))\subset\beta(\bbb{G}_q,v)$.
Since $v$ has a self-arc in $\bbb{G}_q$, it follows that $v\in\beta(\bbb{G}_q,v)$.
Therefore, $\beta(\bbb{G}_q,v)$ is a strictly proper subset of $\scr{V}$ which contains $v$
and is $\beta(\bbb{G}_p,\cdot)$ - invariant. But this is impossible since $\bbb{G}_p$ is sunk at $v$.
\hfill$\qed$

\begin{lemma}
Suppose that $m>1$ and let $\bbb{G}_{p_1},\bbb{G}_{p_2},\ldots,\bbb{G}_{p_k}$ be a finite sequence of
graphs in $\scr{G}_{sa}$ which are all sunk at $v$.
If $k \leq m-1$,
then $v$ has at least $k+1$ neighbors in $\bbb{G}_{p_k}\circ\bbb{G}_{p_{k-1}}\circ\cdots\circ\bbb{G}_{p_1}$.
If $k \ge m-1$,
then the composition $\bbb{G}_{p_k}\circ\bbb{G}_{p_{k-1}}\circ\cdots\circ\bbb{G}_{p_1}$ is strongly sunk at $v$.
\label{kobs}\end{lemma}

\noindent
{\bf Proof of Lemma \ref{kobs}:}
First consider the case when $k\le m-1$.
To prove that $v$ has at least $k+1$ neighbors in
$\bbb{G}_{p_k}\circ\bbb{G}_{p_{k-1}}\circ\cdots\circ\bbb{G}_{p_1}$, suppose the contrary,
namely that $v$ has at most $k$ neighbors in $\bbb{G}_{p_k}\circ\bbb{G}_{p_{k-1}}\circ\cdots\circ\bbb{G}_{p_1}$.
Since each $\bbb{G}_{p_i}$ has self-arcs at all vertices, $i\in\mathbf{k}$,
for each $j\in\{0,1,\ldots,k-1\}$, the arcs in $\bbb{G}_{p_k}\circ\bbb{G}_{p_{k-1}}\circ\cdots\circ\bbb{G}_{p_{k-j}}$
must all be arcs in
$\bbb{G}_{p_k}\circ\bbb{G}_{p_{k-1}}\circ\cdots\circ\bbb{G}_{p_1}$.
Since $k\leq m-1$, it follows that
$\beta(\bbb{G}_{p_k}\circ\bbb{G}_{p_{k-1}}\circ\cdots\circ\bbb{G}_{p_1},v)$ must be a
strictly proper subset of $\scr{V}$, so is
$\beta(\bbb{G}_{p_k}\circ\bbb{G}_{p_{k-1}}\circ\cdots\circ\bbb{G}_{p_{k-j}},v)$
 for all $j\in\{0,1,\ldots,k-1\}$.
By Lemma \ref{func2}, $\beta(\bbb{G}_{p_k}\circ\bbb{G}_{p_{k-1}}\circ\cdots\circ\bbb{G}_{p_{k-j+1}},v)$
is a strictly proper subset of
$\beta(\bbb{G}_{p_k}\circ\bbb{G}_{p_{k-1}}\circ\cdots\circ\bbb{G}_{p_{k-j}},v)$
for  all $j\in\{1,\ldots,k-1\}$.
 In view of this, each containment in
the ascending chain
\begin{eqnarray*}
\beta(\bbb{G}_{p_k},v) &\subset& \beta(\bbb{G}_{p_k}\circ\bbb{G}_{p_{k-1}},v)\subset \cdots \\
&\subset&
\beta(\bbb{G}_{p_k}\circ\bbb{G}_{p_{k-1}}\circ\cdots\circ\bbb{G}_{p_1},v)
\end{eqnarray*}
is strict. Since $\beta(\bbb{G}_{p_k},v)$ has at least two vertices in it,
$\beta(\bbb{G}_{p_k}\circ\bbb{G}_{p_{k-1}}\circ\cdots\circ\bbb{G}_{p_1},v)$ must have at least
$k+1$ vertices in it. But this is impossible because of the hypothesis that
$v$ has at most $k$ neighbors in $\bbb{G}_{p_k}\circ\bbb{G}_{p_{k-1}}\circ\cdots\circ\bbb{G}_{p_1}$.
Therefore, if $\bbb{G}_{p_1},\bbb{G}_{p_2},\ldots,\bbb{G}_{p_k}$ are all sunk at $v$ and $k \leq m-1$,
then $v$ has at least $k+1$ neighbors in $\bbb{G}_{p_k}\circ\bbb{G}_{p_{k-1}}\circ\cdots\circ\bbb{G}_{p_1}$.

Next consider the case when {\color{black}$k\ge m-1$}. In view of the preceding, in the case when $k=m-1$,
the vertex $v$ has at least $m$ neighbors in the composition
$\bbb{G}_{p_{m-1}}\circ\bbb{G}_{p_{m-2}}\circ\cdots\circ\bbb{G}_{p_1}$.
But a vertex in an $m$-vertex graph can have at most $m$
neighbors. Thus, in this case, $v$ has $m$ neighbors, which implies that
$\bbb{G}_{p_{m-1}}\circ\bbb{G}_{p_{m-2}}\circ\cdots\circ\bbb{G}_{p_1}$ is strongly sunk at $v$.
In the case when {\color{black}$k\ge m$}, since each graph considered here has self-arcs at all vertices,
the arcs of $\bbb{G}_{p_{m-1}}\circ\bbb{G}_{p_{m-2}}\circ\cdots\circ\bbb{G}_{p_1}$ must all be arcs in
$\bbb{G}_{p_k}\circ\bbb{G}_{p_{k-1}}\circ\cdots\circ\bbb{G}_{p_1}$, so
$v$ has $m$ neighbors in $\bbb{G}_{p_k}\circ\bbb{G}_{p_{k-1}}\circ\cdots\circ\bbb{G}_{p_1}$
for all {\color{black}$k\ge m$}. Thus, the composition
$\bbb{G}_{p_{k}}\circ\bbb{G}_{p_{k-1}}\circ\cdots\circ\bbb{G}_{p_1}$ is strongly sunk at $v$
for all {\color{black}$k\ge m-1$}.
\hfill$\qed$


To prove Proposition \ref{p}, we will also make use of the following idea.
By a {\em route} over a given sequence of graphs $\mathbb{G}_1,\mathbb{G}_2,\ldots,\mathbb{G}_q$  in $\scr{G}_{sa}$
 is meant a sequence of vertices $i_0, i_1,\ldots, i_q$ such that
  $(i_{k-1},i_k)$  is an arc in $\mathbb{G}_k$ for all $k\in\mathbf{q}$.
  A route over a sequence of graphs which are all the same graph $\mathbb{G}$, is thus a
  walk in $\mathbb{G}$.

  The definition of a route implies that if $i_0, i_1,\ldots, i_q$ is a route over
$\mathbb{G}_1,\mathbb{G}_2,\ldots,\mathbb{G}_q$ and  $i_q,i_{q+1},\ldots,$ $i_p$ is a route over
$\mathbb{G}_q,\mathbb{G}_{q+1},\ldots,\mathbb{G}_p$,
then  the `concatenated' sequence $i_0, i_1,\ldots,i_{q-1}, i_q,i_{q+1},\ldots, i_p$
is a route over $\mathbb{G}_1,\mathbb{G}_2,\ldots,\mathbb{G}_{q-1},\mathbb{G}_{q}, \mathbb{G}_{q+1},\ldots,\mathbb{G}_p$.
This remains true if more than two sequences are concatenated.

Note that the definition of composition in $\scr{G}_{sa}$  implies that if $i=i_0, i_1,\ldots, i_q=j$ is
a route over a sequence $\mathbb{G}_1,\mathbb{G}_2,\ldots,\mathbb{G}_q$,
then $(i,j)$ must be an arc in the composed graph $\mathbb{G}_q\circ\mathbb{G}_{q-1}\circ\cdots \circ\mathbb{G}_1$.
The definition of composition also implies the converse, namely that  if $(i,j)$ is an arc in $\mathbb{G}_q\circ\mathbb{G}_{q-1}\circ\cdots \circ\mathbb{G}_1$,
then there must exist vertices $i_1,\ldots, i_{q-1}$ for which $i=i_0, i_1,\ldots, i_q=j$ is a route over
$\mathbb{G}_1,\mathbb{G}_2,\ldots,\mathbb{G}_q$.

Suppose that $\bbb{G}_{\tau_1},\bbb{G}_{\tau_2},\ldots,\bbb{G}_{\tau_p}$ is a subsequence of
$\mathbb{G}_1,\mathbb{G}_2,\ldots,\mathbb{G}_q$ in $\scr{G}_{sa}$, $p\le q$, and that
$i_0, i_1,\ldots, i_p$ is a route over
$\bbb{G}_{\tau_1},\bbb{G}_{\tau_2},\ldots,\bbb{G}_{\tau_p}$. Then,
since each $\bbb{G}_i$, $i\in\mathbf{q}$, has self-arcs at all vertices,
there must exist a route over $\mathbb{G}_1,\mathbb{G}_2,\ldots,\mathbb{G}_q$ which
contains $i_0, i_1,\ldots, i_p$ as a subsequence.

More can be said.

\begin{lemma}
{\color{black}{\rm (Lemma 4 in \cite{tac})}}
Let $S_1,S_2,\ldots S_q $ be a sequence of $m\times m$  stochastic matrices
with  graphs
$\mathbb{G}_1,\;\mathbb{G}_2,\ldots,\mathbb{G}_q$  in $\scr{G}_{sa}$ respectively. If  $j=i_0, i_1,\ldots, i_q=i$
 is a
route over the sequence
$\mathbb{G}_1,\;\mathbb{G}_2,\ldots,\mathbb{G}_q$, then the   matrix product $P_{i_q}P_{i_{q-1}}\cdots P_{i_0}$
is a component of the   $ij$th block entry of
$$M=P(S_q\otimes I)P(S_{q-1}\otimes I)\cdots  P(S_1\otimes I)P$$ \label{gum}\end{lemma}

{\color{black}
To proceed, we call matrices of the form
$$\mu(P_1,P_2,\ldots ,P_m) = \sum_{i=1}^{d}\lambda_i P_{h_i(1)}P_{h_i(2)}\cdots P_{h_i(q_i)}$$
{\em projection matrix polynomials}
where $q_i$ and $d$ are positive integers, $\lambda_i $ is a real positive
   number, and for each $j\in\{1,2,\ldots,q_i\}$,
 $h_i(j)$  is an integer in $\mathbf{m}$.
 A  projection matrix polynomial
 $\mu(P_1,P_2,\ldots ,P_m)$  is {\em complete}  if it has a    component
   $P_{h_i(1)}P_{h_i(2)}\cdots P_{h_i(q_i)}$ such that for some fully populated set $\scr{E}$,
each of the    projection matrices $P_j$, $j\in\scr{E}$,
appears in the component at least once.
}

We are now in a position to prove Proposition \ref{p} and
the sufficiency of Theorem \ref{main1}.

{\color{black}
In the following proof, we will make use of the fact that
for any $m\times m$ block matrix $Q = \matt{Q_{ij}}$
whose $ij$th entry is a projection matrix polynomial,
if \rep{assmp} holds and at least one entry in each block
row of $Q$ is complete, then $Q$ is a contraction in the mixed
matrix norm (see Proposition 1 in \cite{tac}).
}

\noindent
{\bf Proof of Proposition \ref{p}:}
Let $v$ be any vertex in $\scr{V}=\mathbf{m}$ and
$\mathbb{G}_i =\gamma(S_{(i+1)l})\circ\gamma(S_{(i+1)l-1})\circ\cdots\circ\gamma(S_{il+1})$,
$i\in\{0,1,2,\ldots\}$.
Since the sequence $\gamma(S_1),\gamma(S_2),\ldots$ is repeatedly jointly $\scr{D}$-connected
by subsequences of length $l$,
each $\mathbb{G}_i$ is $\scr{D}$-connected.
Let $\scr{V}_i(v) = \{v\}\cup \beta(\bbb{G}_i,v)\cup\cdots\cup\beta^{m-1}(\bbb{G}_i,v)$
be the vertex set of $\mathbb{G}_i$'s
sunk graph generated by $v$. It follows from Lemma \ref{easy} that $\scr{V}_i(v)$ must be fully populated.
Since $\scr{V}_i(v)\subset \scr{V}$ and $\scr{V}$ is a finite set,
the total number of possibly distinct $\scr{V}_i(v)$ is finite.
We use $h(v)$ to denote this finite number.

Let $d_i$ denote the cardinality of $\scr{V}_i(v)$ and let $\bar d$ be the maximum of all possible $d_i$.
Since $\scr{V}_i(v)\subset \scr{V}$, it follows that $\bar d\leq m$.
Set $r=\frac{1}{2}(\bar d-1)\bar d$ and $\tau(v)=(r-1)h(v)+1$. Then, there must exist one set
which appears in the sequence $\scr{V}_1(v),\scr{V}_2(v),\ldots,\scr{V}_\tau(v)$ at least $r$ times.
We use $\scr{U}$ to denote this set and let $d$ be its cardinality.
Let $\bbb{G}_{\tau_1},\bbb{G}_{\tau_2},\ldots,\bbb{G}_{\tau_q}$ be the subsequence of
$\bbb{G}_1,\bbb{G}_2,\ldots,\bbb{G}_{\tau(v)}$
which includes every graph in the sequence $\bbb{G}_1,\bbb{G}_2,\ldots,\bbb{G}_{\tau(v)}$
whose vertex set of the sunk graph generated by $v$
is $\scr{U}$.
Then, there holds $\tau(v)\ge q\geq r \geq \frac{1}{2}(d-1)d$.

We will show that there exists at least one complete 
block entry in the $v$th row of
$M=P(S_p\otimes I)P(S_{p-1}\otimes I)\cdots P(S_1\otimes I)P$ for any $p\ge \tau(v)$.
In the case when $d=1$, it must be true that $\scr{U}=\{v\}$. By Lemma \ref{easy}, there holds $\scr{P}_v=\{0\}$.
Then, it is straightforward to verify that each block entry of $v$th row of $M$ is complete.
Next we consider the case when $d>1$.

Let $\bbb{H}_1, \bbb{H}_2,\ldots,\bbb{H}_q$ be a sequence of
graphs such that $\bbb{H}_1=\bbb{G}_{\tau_q}, \bbb{H}_2=\bbb{G}_{\tau_{q-1}}, \ldots, \bbb{H}_q=\bbb{G}_{\tau_1}$.
It follows immediately that
if $i_q,i_{q-1},\ldots,i_0$ is a route over $\bbb{H}_q, \bbb{H}_{q-1},\ldots,\bbb{H}_1$,
then $i_q,i_{q-1},\ldots,i_0$ is a route over $\bbb{G}_{\tau_1}, \bbb{G}_{\tau_2},\ldots,\bbb{G}_{\tau_q}$.
Set $\Sigma_k= 1+2+\cdots + k = \frac{1}{2}k(k+1)$ for each $k\in\{1,2,\ldots,d-1\}$.
Partition the sequence $\bbb{H}_1, \bbb{H}_2,\ldots,\bbb{H}_q$ into
$d-1$ successive subsequences:
$\scr{H}_1= \bbb{H}_1$,
$\scr{H}_2= \bbb{H}_2,\bbb{H}_3$,
$\scr{H}_3= \bbb{H}_4,\bbb{H}_5,\bbb{H}_6$,
$\ldots,$
$\scr{H}_{d-2} = \bbb{H}_{\Sigma_{d-3}+1},\bbb{H}_{\Sigma_{d-3}+2},\ldots,\bbb{H}_{\Sigma_{d-2}}$,
and $\scr{H}_{d-1} = \bbb{H}_{\Sigma_{d-2}+1},\bbb{H}_{\Sigma_{d-2}+2},\ldots,\bbb{H}_q$.
Note that each $\scr{H}_i$ is of length $i$ except for the last which must be of length
$q-\Sigma_{d-2}\geq d-1$.

Let $\scr{U}=\{i_1,i_2,\ldots,i_d\}$ and $i_1=v$.
Since $\scr{U}$ is the vertex set of the sunk graph generated by $v$ for each $\bbb{H}_k$, $k\in\mathbf{q}$,
and $d>1$, it follows that
$i_1$ must have a neighbor $i_2$ in $\bbb{H}_1$
which is not itself. Then, there is a route over $\scr{H}_1$ from $i_2$ to $i_1$.
Note that $d\le m$.
By Lemma \ref{kobs}, it follows that $i_2$ has at least three neighbors
in $\bbb{H}_{2}\circ \bbb{H}_3$.
Thus, $i_2$ must have a neighbor $i_3$ in $\bbb{H}_{2}\circ \bbb{H}_3$
which is not in the set $\{i_1,i_2\}$.
Then, there is a route over $\scr{H}_2$ from $i_3$ to $i_2$.
By repeating this argument for all $\scr{H}_k$,
one could get a sequence of distinct vertices $i_1,i_2,\ldots,i_d$ such that
there is a route over $\scr{H}_k$ from $i_{k+1}$ to $i_{k}$ for all $k\in\{1,2,\ldots,d-1\}$.

For each $\scr{H}_k = \bbb{H}_{\Sigma_{k-1}+1},\ldots,\bbb{H}_{\Sigma_{k}}$,
$k\in\{1,2,\ldots,d-2\}$, we use
$i_{k+1}=j_{\Sigma_{k}},j_{\Sigma_{k}-1},\ldots,j_{\Sigma_{k-1}}=i_k$
to denote its route from $i_{k+1}$ to $i_{k}$.
For $\scr{H}_{d-1} = \bbb{H}_{\Sigma_{d-2}+1},\ldots,\bbb{H}_q$,
we use $i_{d}=j_{q},j_{q-1},\ldots,j_{\Sigma_{d-2}}=i_{d-1}$
to denote its route from $i_{d}$ to $i_{d-1}$.
Thus, $j_q,j_{q-1},\ldots,j_0$ must be a route
over the overall sequence $\bbb{H}_q,\bbb{H}_{q-1},\ldots,\bbb{H}_1$. In particular,
$j_{\Sigma_{k-1}}=i_k$ for all $k\in\{1,2,\ldots,d-1\}$ and  $j_q=i_d$.

From the preceding, it follows that
$i_d=j_q,j_{q-1},\ldots,j_0 = i_1=v$  is a route over $\bbb{G}_{\tau_1},\bbb{G}_{\tau_2},\ldots,\bbb{G}_{\tau_q}$.
Since $\bbb{G}_{\tau_1},\bbb{G}_{\tau_2},\ldots,\bbb{G}_{\tau_q}$ is a subsequence
of $\bbb{G}_1,\bbb{G}_2,\ldots,\bbb{G}_p$ and
all the graphs considered here have self-arcs at all vertices,
there must be a route $k_0,k_1,\ldots,k_s=v$ over $\bbb{G}_1, \bbb{G}_2,\ldots,\bbb{G}_p$
which contains $i_d=j_q,j_{q-1},\ldots,j_0 = i_1 =v$ as a subsequence.
In view of Lemma \ref{gum}, the matrix product $P_{k_{0}}P_{k_{1}}\cdots P_{k_{s}}$
must be a component of the $(vk_0)$th block entry of
$$M=P(S_p\otimes I)P(S_{p-1}\otimes I)\cdots  P(S_1\otimes I)P$$
But $\scr{U}=\{i_1,i_2,\ldots,i_d\}$ is a fully populated subset
and each element of $\scr{U}$ appears in $k_0,k_{1},\ldots,k_s$  at least once,
so the $(vk_0)$th block entry of $M$ must be complete.
Since this reasoning applies for any $v\in\mathbf{m}$,
at least one entry in each row of $M$ is a complete projection matrix polynomial if $p\ge \tau(v)$ for all $v\in\mathbf{m}$.
It follows from Proposition 1 of \cite{cdc13.1} that $M$ is a contraction in the mixed matrix norm.
\hfill$\qed$

\noindent
{\bf Proof of Theorem \ref{main1} (Sufficiency):}
Since the sequence $\bbb{N}(1),\bbb{N}(2),\ldots$ is repeatedly jointly $\scr{D}$-connected
and $\bbb{N}(t)=\gamma(F(t))$ for all $t\geq 1$, the sequence $\gamma(F(1)),\gamma(F(2)),\ldots$
is repeatedly jointly $\scr{D}$-connected.
Without loss of generality, suppose that
the sequence $\gamma(F(1)),\gamma(F(2)),\ldots$ is repeatedly jointly $\scr{D}$-connected
by subsequences of length $l$.
By Proposition \ref{p}, there is a finite positive integer $\tau$, not depending on $P$, such that
for any  $p\geq \tau$,  the matrix
$P(F(p)\otimes I)P(F(p-1)\otimes I)\cdots  P(F(1)\otimes I)P$
  is a contraction in the mixed  matrix norm.
Let $q$ be the smallest integer such that $ql\geq \tau$. Set $\omega=ql$.
Then, for each $k\geq 0$, the sequence $\gamma(F(k\omega+1)),\gamma(F(k\omega+2)),\ldots$
is repeatedly jointly $\scr{D}$-connected by subsequence of length $l$ and thus
the matrix
$P(F((k+1)\omega)\otimes I)P(F((k+1)\omega-1)\otimes I)\cdots  P(F(k\omega+1)\otimes I)P$
  is a contraction in the mixed  matrix norm.
Since directed graphs
 in $\scr{G}_{sa}$ are bijectively related to flocking matrices,
 the set $\scr{F}_\omega $ of distinct subsequences
 $F(k\omega +1),F(k\omega+2),\ldots, F((k+1)\omega)$, $k\geq 0$,
 encountered along any trajectory of \rep{sys}, must be a finite
set and thus compact.
Thus,
\begin{eqnarray*}
\lambda &=& \Bigg(\sup_{k\in\{0,1,2,\ldots\}} \big\|P(F((k+1)\omega)\otimes I)\cdots \\
&&  P(F(k\omega+2)\otimes I)P(F(k\omega+1)\otimes I)P\big\|\Bigg)^{\frac{1}{\omega}} <1
\end{eqnarray*}
This thus completes the proof of Theorem \ref{main1}.
\hfill$\qed$

\subsection{Nonuniqueness} \label{nonunique}

We now turn to the case when $Ax=b$ has multiple solutions (Theorem \ref{main2}).
To deal with this case we will ``quotient out'' the subspace $\bigcap_{i=1}^m\scr{P}_i$,
much as was done in \cite{cdc13.1}, thereby decomposing the problem into two
parts - one to which the result for the uniqueness case (Theorem \ref{main1}) is applicable
and the other to which a result for standard consensus is applicable.
The steps involved in doing this make use of the following lemma.

\begin{lemma}
Let $Q'$ be any matrix whose columns form an orthonormal basis for the orthogonal complement of the
subspace $\bigcap_{i=1}^m\scr{P}_i$ and define $\bar{P}_i = QP_iQ'$, $i\in\mathbf{m}$.
Then, the following statements are true.
\mbox{}

1. Each $\bar{P}_i$, $i\in\mathbf{m}$, is an orthogonal projection matrix.

2. Each $\bar{P}_i$, $i\in\mathbf{m}$,  satisfies $QP_i = \bar{P}_i Q$.


 3. For any nonempty subset $\scr{E}\subset\mathbf{m}$, there holds $\bigcap_{i\in\scr{E}} \bar{\scr{P}}_i =\{0\}$ if and only if $\scr{E}$ is fully populated.


\label{vlad}\end{lemma}

\noindent{\bf Proof of Lemma \ref{vlad}:}
Note that $\bar{P}_i^2 =QP_iQ'QP_iQ'=QP_i^2Q' =QP_iQ' = \bar{P}_i$, $i\in\mathbf{m}$, so
each $\bar{P}_i$ is idempotent; since each $\bar{P}_i$ is symmetric, each must be an orthogonal projection matrix.
Thus, property 1 is true.

Since $\ker Q = \bigcap_{i=1}^m \scr{P}_i$, it must be true that  $\ker Q \subset  \scr{P}_i$,
$i\in\mathbf{m}$. Thus,
$P_i\ker Q =\ker Q$, $i\in\mathbf{m}$.
Therefore,  $QP_i\ker Q = 0$, so $\ker Q\subset \ker QP_i$.  This and the fact that $Q$ has
 linearly independent rows mean that the equation $QP_i = XQ$ has a unique solution $X$. It follows that $X = QP_iQ'$,
 so $X=\bar{P}_i$.  Therefore, property 2 is true.


Suppose first that $\scr{E}\subset\mathbf{m}$ is fully populated.
Then, $\bigcap_{i\in\scr{E}}\scr{P}_i=\bigcap_{i=1}^m\scr{P}_i$.
Pick $x\in\bigcap_{i\in\scr{E}}\bar\scr{P}_i$.  Then, $x\in\bar{\scr{P}}_i$, $i\in\scr{E}$, so
 there exist $w_i$ such that $x=\bar{P}_iw_i$, $i\in\scr{E}$.  Set   $y=Q'x$ in which case $x=Qy$; thus,  $y
  = Q'\bar{P}_iw_i$, $i\in\scr{E}$. In view of property 2 of Lemma \ref{vlad},
 $y= P_iQ'w_i$, $i\in\scr{E}$,  so $y\in\bigcap_{i\in\scr{E}}\scr{P}_i$.  Thus, $Qy = 0$.  But $x=Qy$, so $x=0$.
 Therefore, $\bigcap_{i\in\scr{E}} \bar{\scr{P}}_i =\{0\}$.
Suppose next that $\bigcap_{i\in\scr{E}} \bar{\scr{P}}_i =\{0\}$.
 Pick $y\in\bigcap_{i\in\scr{E}}\scr{P}_i$.
  Then, $y\in\scr{P}_i$, $i\in\scr{E}$, so
 there exist $v_i$ such that $y=P_iv_i$,  $i\in\scr{E}$.
 Set  $x=Qy$  in which case $y=Q'x$; thus,  $x  = QP_iv_i$, $i\in\scr{E}$.
  In view of property 2 of Lemma \ref{vlad},
 $x= \bar P_iQv_i$, $i\in\scr{E}$,  so $x\in\bigcap_{i\in\scr{E}}\bar\scr{P}_i$ and $x=0$.
 But $x=Qy$, so $y\in\ker Q = \bigcap_{i=1}^m \scr{P}_i$.
 This implies that $\bigcap_{i\in\scr{E}}\scr{P}_i=\bigcap_{i=1}^m\scr{P}_i$
 and thus $\scr{E}$ is fully populated.
Therefore, property 3 is true.
\hfill$\qed$

We are now in a position to prove Theorem \ref{main2}.

\noindent
{\bf Proof of Theorem \ref{main2}:}
Property 2 of Lemma \ref{vlad} implies that $QP_iP_j = \bar{P}_i\bar{P}_jQ$ for all
$i,j\in\mathbf{m}$.  If
we define $\bar{y}_i = Qy_i$, $i\in\mathbf{m}$,   then from \rep{na1},
\eq{\bar{y}_i(t+1) =
\frac{1}{m_i(t)}\bar{P}_i\sum_{j\in\scr{N}_i(t)}\bar{P}_j\bar{y}_j(t),\;\;\;t\geq 1\label{nna1}}
Define $z_i = y_i-Q'\bar{y}_i$, $i\in\mathbf{m}$.
Note that $Qz_i = Qy_i - \bar{y}_i$, so $Qz_i = 0$, $i\in\mathbf{m}$.  Thus, $z_i(t)\in\bigcap_{j=1}^m \scr{P}_j$,
$i\in\mathbf{m}$.
It follows that $P_jz_i(t) = z_i(t)$, $i,j\in\mathbf{m}$.
Moreover,  from property 2 of Lemma \ref{vlad}, $P_iQ' = Q'\bar{P}_i$.  These expressions and \rep{nna1}
imply that
\eq{z_i(t+1) =
\frac{1}{m_i(t)}\sum_{j\in\scr{N}_i(t)}z_j(t),\;\;\;t\geq 1\label{dinner}}
These equations are the update equations for the standard consensus problem treated in \cite{reachingp1} and elsewhere for the case
 when $z_i$'s are scalars.

Since by assumption, $A\neq 0$, the matrix $Q$ defined in the statement of Lemma \ref{vlad} is not the
zero matrix and so the subsystem defined by \rep{nna1} has positive dimension.
Note that \rep{nna1} has exactly the same form as \rep{na1} except for the  $\bar{P}_i$'s which replace the $P_i$'s.
But in view of Lemma \ref{vlad},  the $\bar{P}_i$'s are also orthogonal  projection matrices  and
$\bigcap_{i=1}^m\bar{\scr{P}}_i = \{0\}$.
Therefore, all $x_i(t)$ defined by  \rep{a1} converge to the same limit exponentially fast
if and only if all $\bar y_i(t)$ defined by  \rep{nna1} converge to $0$ exponentially fast
and all $z_i(t)$ defined by  \rep{dinner} converge to the same limit $z^*\in\bigcap_{i = 1}^m\scr{P}_i$ exponentially fast.
The limit of $x_i(t)$ is $x^* + z^*$ which solves $Ax=b$.

Since $\bigcap_{i=1}^m\bar{\scr{P}}_i = \{0\}$, Theorem \ref{main1}
  is applicable to the system of iterations \rep{nna1}.
Property 3 of Lemma \ref{vlad} implies that a subset
$\scr{E}\subset\mathbf{m}$ is fully populated with respect to $\{P_1,P_2,\ldots,P_m\}$
(i.e., $\bigcap_{i\in\scr{E}}\scr{P}_i = \bigcap_{i=1}^m\scr{P}_i$)
if and only if $\scr{E}$ is fully populated with respect to $\{\bar P_1,\bar P_2,\ldots,\bar P_m\}$
(i.e., $\bigcap_{i\in\scr{E}}\bar\scr{P}_i = \bigcap_{i=1}^m\bar\scr{P}_i$).
Property 3 of Lemma \ref{vlad} implies that a subset
$\scr{E}\subset\mathbf{m}$ is partially populated with respect to $\{P_1,P_2,\ldots,P_m\}$
(i.e., $\bigcap_{i=1}^m\scr{P}_i$ is a proper subset of $\bigcap_{i\in\scr{E}}\scr{P}_i$)
if and only if $\scr{E}$ is partially populated with respect to $\{\bar P_1,\bar P_2,\ldots,\bar P_m\}$
(i.e., $\bigcap_{i=1}^m\bar\scr{P}_i$ is a proper subset of $\bigcap_{i\in\scr{E}}\bar\scr{P}_i$).
Therefore, a directed graph $\bbb{G}\in\scr{G}_{sa}$ with vertex set $\scr{V}=\mathbf{m}$
is $\scr{D}$-connected if and only if it is $\bar\scr{D}$-connected,
where $\bar\scr{D}$ is the collection of all partially populated subsets of $\scr{V}$ with respect to
$\{\bar P_1,\bar P_2,\ldots,\bar P_m\}$.
It then follows that $\bar y_i(t)\rightarrow 0$ exponentially fast as $t\rightarrow\infty$
if and only if
the sequence of neighbor graphs $\bbb{N}(1),\bbb{N}(2),\ldots$ is repeatedly jointly $\scr{D}$-connected.

It is known that for the scalar case, a necessary and sufficient
 condition for  all $z_i(t)$ defined by \rep{dinner}
 to converge exponentially fast to the same value is that
 the sequence  $\bbb{N}(1),\bbb{N}(2),\ldots$ is repeatedly jointly rooted
 \cite{luc,cdc14}.
But since   the  vector update  \rep{dinner} decouples into
 $n$ independent scalar update equations,  the convergence conditions for the scalar equations apply  without any change
 to the vector case  as well. Thus,  all $z_i(t)$ converge  exponentially  fast
  to the same limit  $z^*\in\bigcap_{i = 1}^m\scr{P}_i$ if and only if
  $\bbb{N}(1),\bbb{N}(2),\ldots$ is repeatedly jointly rooted.
It follows from the preceding discussion that for the case when $Ax=b$ has multiple solutions and $A\neq 0$,
all $x_i(t)$ defined by  \rep{a1} converge to the same solution to $Ax=b$
if and only if
the sequence of neighbor graphs $\bbb{N}(1),\bbb{N}(2),\ldots$ is repeatedly jointly $\scr{D}$-connected
and is repeatedly jointly rooted.
\hfill$\qed$

\section{Discussions} \label{discuss}

In this section, we present some additional results.
A convergability issue is addressed in \S \ref{compare},
a further result on the notion of $\scr{D}$-connectivity is given in \S \ref{equal},
{\color{black}
a lower bound on the convergence rate is derived in \S \ref{rate},
and the case when $Ax=b$ does not have a solution is discussed in \S \ref{least}.
}

\subsection{Comparison with Consensus} \label{compare}


In \cite{pieee},
 a set $\scr{S}$ of $m\times m$ stochastic matrices   with graphs in $\scr{G}_{sa}$ is called
{\em
convergable} if for every compact subset $\scr{C}\subset \scr{S}$ and every infinite sequence
$S_1,S_2,\ldots$  of matrices $S_i$  from  $\scr{C}$, the matrix product
$S_tS_{t-1}\cdots S_1$
converges as $t\rightarrow\infty$ to  a rank one  matrix of the form $\mathbf{1}c$ where $\mathbf{1}$
is a column vector in $\R^n$ whose entries all equal $1$ and $c$ is a row vector in $\R^{1\times n}$.
It is known that $\scr{S}$ is convergable if and only if  every matrix in $\scr{S}$ has a rooted graph
(see Theorem 4 of \cite{pieee}).
It is natural to wonder what the analogous condition would be for the types of matrix products under consideration here.
With this in mind, let us say that a set $\scr{S}$ of $m\times m$ stochastic matrices
with graphs in $\scr{G}_{sa}$ is
{\em $P$-convergable} if for every compact subset $\scr{C}\subset \scr{S}$ and every infinite sequence
$S_1,S_2,\ldots$  of matrices $S_i$  from  $\scr{C}$, the matrix
$P(S_t\otimes I)P(S_{t-1}\otimes I)\cdots  P(S_1\otimes I)P$
converges to the zero matrix as $t\rightarrow\infty$.

\begin{corollary}
Suppose that \rep{assmp} holds.
 A set $\scr{S}$  of stochastic matrices with graphs in $\scr{G}_{sa}$ is   $P$-convergable if and only of every matrix
 in $\scr{S}$ has a $\scr{D}$-connected graph.
\label{pcon}\end{corollary}
The necessity is an immediate consequence of Proposition \ref{ns}
and the sufficiency is a direct consequence of Theorem \ref{main1}.
The corollary implies that in the case when \rep{assmp} holds,
the set $\scr{A}$ consisting of all stochastic matrices  with graphs in $\scr{G}_{sa}$ which is $\scr{D}$-connected is
the largest subset of $m\times m$ stochastic matrices with positive diagonal entries
whose compact subsets are all $P$-convergable.
The set $\scr{A}$ itself is not $P$-convergable because it is not closed and thus not compact.

\subsection{Comparison with Strong Connectivity} \label{equal}

As noted in \S \ref{pconnect}, the notion of $\scr{D}$-connectivity is less restrictive than strong connectivity.
We present here a condition under which the two properties are equivalent.

\begin{lemma}
Suppose that any proper subset of $\scr{V}$ is partially populated.
Then, a graph $\bbb{G}$ in $\scr{G}_{sa}$ is $\scr{D}$-connected if and only if it is strongly connected.
\label{ps}\end{lemma}

\noindent
{\bf Proof of Lemma \ref{ps}:}
Since each strongly connected graph must be $\scr{D}$-connected, it is enough to prove the sufficiency.
Suppose that $\bbb{G}$ in $\scr{G}_{sa}$ is a graph which is $\scr{D}$-connected and
let $v$ be any vertex in $\scr{V}=\mathbf{m}$.
Since any proper subset of $\scr{V}$ is partially populated,
it follows that $\scr{V}(v) = \{v\}\cup\beta(\bbb{G},v)\cup \cdots \cup \beta^{m-1}(\bbb{G},v)$ must equal $\scr{V}$.
This implies that the sunk graph generated by $v$ is $\bbb{G}$ itself. Thus, $v$ is reachable from
each other vertex of $\bbb{G}$.
Since this reasoning applies for any $v\in\scr{V}$, it follows that $\bbb{G}$ must be strongly connected.
\hfill$\qed$

{\color{black}

\subsection{An Improved Bound on Convergence Rate} \label{rate}

In this section, we will present an improved lower bound on convergence of the algorithm \rep{a1} in the case
when 
the neighbor graph sequence is repeatedly jointly strongly connected,
which is tighter than that in \cite{tac}.


\begin{proposition}
For any set of $\frac{1}{2}(m-1)m$ or more $m\times m$ stochastic matrices
$S_1,S_2,\ldots,S_q$ whose graphs are all strongly connected in $\scr{G}_{sa}$, the matrix
$P(S_q\otimes I)P(S_{q-1}\otimes I)\cdots  P(S_1\otimes I)P$ is a contraction
in the mixed matrix norm.
\label{lowerxxx}\end{proposition}

{\bf Proof of Proposition \ref{lowerxxx}:}
Set $\mathbb{G}_i =\gamma(S_i),\;i\in\{1,2,\ldots,p\}$.
Let $\bbb{H}_1, \bbb{H}_2,\ldots,\bbb{H}_q$ be a sequence of
graphs such that $\bbb{H}_1=\bbb{G}_q, \bbb{H}_2=\bbb{G}_{q-1}, \ldots, \bbb{H}_q=\bbb{G}_1$.
It should be clear that each $\bbb{H}_i$ is a strongly connected graph in $\scr{G}_{sa}$.
It should also be clear that if $i_q,i_{q-1},\ldots,i_0$ is a route over $\bbb{H}_q, \bbb{H}_{q-1},\ldots,\bbb{H}_1$,
then $i_q,i_{q-1},\ldots,i_0$ is a route over $\bbb{G}_1, \bbb{G}_2,\ldots,\bbb{G}_q$.

Set $\Sigma_k= 1+2+\cdots + k = \frac{1}{2}k(k+1)$ for $k\in\{1,2,\ldots,m-1\}$.
Partition the sequence $\bbb{H}_1, \bbb{H}_2,\ldots,\bbb{H}_q$ into
$m-1$ successive subsequences
$\scr{H}_1= \bbb{H}_1$,
$\scr{H}_2= \bbb{H}_2,\bbb{H}_3$,
$\scr{H}_3= \bbb{H}_4,\bbb{H}_5,\bbb{H}_6$,
$\ldots$ ,
$\scr{H}_{m-2} = \bbb{H}_{\Sigma_{m-3}+1},\bbb{H}_{\Sigma_{m-3}+2},\ldots,\bbb{H}_{\Sigma_{m-2}}$,
$\scr{H}_{m-1} = \bbb{H}_{\Sigma_{m-2}+1},\bbb{H}_{\Sigma_{m-2}+2},\ldots,\bbb{H}_q$,
each $\scr{H}_i$ of length $i$ except for the last which must be of length
$q-\Sigma_{m-2}\geq m-1$.

Let $i_1$ be any vertex in $\scr{V}=\{1,2,\ldots,m\}$.
Since $\bbb{H}_1$ is strongly connected, $i_1$ must have a neighbor $i_2$ in $\bbb{H}_1$
which is not itself. Then, there is a route over $\scr{H}_1$ from $i_2$ to $i_1$.
By Lemma \ref{kobs} and the fact that strongly connected graphs are sunk at every vertex,
$i_2$ has at least three neighbors
in $\bbb{H}_{2}\circ \bbb{H}_3$.
Thus, $i_2$ must have a neighbor $i_3$ in $\bbb{H}_{2}\circ \bbb{H}_3$
which is not in $\{i_1,i_2\}$.
Then there is a route over $\scr{H}_2$ from $i_3$ to $i_2$.
By repeating this argument for all $\scr{H}_k$,
one could get a sequence of distinct vertices $i_1,i_2,\ldots,i_m$ such that
there is a route over $\scr{H}_k$ from $i_{k+1}$ to $i_{k}$ for all $k\in\{1,2,\ldots,m-1\}$.

For each $\scr{H}_k = \bbb{H}_{\Sigma_{k-1}+1},\ldots,\bbb{H}_{\Sigma_{k}}$,
$k\in\{1,2,\ldots,m-2\}$, we write
$i_{k+1}=j_{\Sigma_{k}},j_{\Sigma_{k}-1},\ldots,j_{\Sigma_{k-1}}=i_k$
for its route from $i_{k+1}$ to $i_{k}$.
For $\scr{H}_{m-1} = \bbb{H}_{\Sigma_{m-2}+1},\ldots,\bbb{H}_q$,
we write $i_{m}=j_{q},j_{q-1},\ldots,j_{\Sigma_{m-2}}=i_{m-1}$
for its route from $i_{m}$ to $i_{m-1}$.
Thus, $j_q,j_{q-1},\ldots,j_0$ must be a route
over the overall sequence $\bbb{H}_q,\bbb{H}_{q-1},\ldots,\bbb{H}_1$. In particular,
$j_{\Sigma_{k-1}}=i_k$ for all $k\in\{1,2,\ldots,m-1\}$ and  $j_q=i_m$.

From the preceding arguments, it should be clear that
$i_m=j_q,j_{q-1},\ldots,j_0 = i_1$  is a route over $\bbb{G}_1, \bbb{G}_2,\ldots,\bbb{G}_q$.
In view of Lemma \ref{gum}, the matrix product $P_{j_{0}}P_{j_{1}}\cdots P_{j_{q}}$
must be a component of the $(i_1i_m)$th block entry of
$$M=P(S_q\otimes I)P(S_{q-1}\otimes I)\cdots  P(S_1\otimes I)P$$
But $i_1,i_2,\ldots,i_m$ are distinct integers
and each of them appears in $j_q,j_{q-1},\ldots,j_0$  at least once,
so the $(i_1i_m)$th block entry of $M$ must be complete.
Since this reasoning applies for any $i_1\in\{1,2,\ldots,m\}$,
at least one block entry in each row of $M$ is a complete projection matrix polynomial.
It then follows that $M$ is a contraction.
\hfill$\qed$

Let $l$ be a positive integer.
A compact subset of $m\times m$ stochastic matrices
with graphs in $\scr{G}_{sa}$ is {\em $l$-compact}
if  the set  $\scr{C}_l$ consisting
 of all sequences
$S_1,S_2,\ldots, S_l$, $S_i\in\scr{C}$, for which $\gamma(S_lS_{l-1}\cdots S_1)$ is strongly connected,
 is nonempty and compact. Thus any nonempty compact subset of  $m\times m$
 stochastic matrices with strongly connected graphs  in $\scr{G}_{sa}$ is
 $1$-compact.

Using Proposition \ref{lowerxxx} and arguments similar to those in \cite{tac}, we can prove the following
result.

\begin{theorem}Suppose that \rep{assmp} holds. Let $l$ be a positive integer.
 Let $\scr{C}$ be an $l$-compact subset of  $m\times m$ stochastic matrices,
and define
\begin{eqnarray*}
    \nonumber \lambda = \Bigg(\sup_{\scr{H}_{\omega}\in\scr{C}_l}\;\sup_{\scr{H}_{\omega-1}\in\scr{C}_l}
\;\cdots \;\sup_{\scr{H}_{1}\in\scr{C}_{l}}\ \ \ \ \ \ \ \ \ \ \ \ \ \ \ \ \ \ \ \ \ \ \ \ \ \ \\
      \ \ \ \ \ \ \ \ \ ||P(Q_{\omega l}\otimes I)P(Q_{\omega l-1}\otimes
 I)\cdots  P(Q_1\otimes I)P|| \Bigg)^{\frac{1}{\omega l}}
\end{eqnarray*}
where $\omega = \frac{1}{2}(m-1)m$ and for  $i\in\{1,2,\ldots, \omega\}$,
and $\scr{H}_i$  is the subsequence $
 Q_{(i-1)l +1},Q_{(i-1)l +2},\ldots, Q_{il}$.   Then, $\lambda <1$, and
 for any infinite sequence
of stochastic matrices $S_1,S_2,\ldots $ in $\scr{C}$ whose graphs form a
sequence $\gamma(S_1),\gamma(S_2),\ldots $ which is repeatedly
jointly strongly connected by contiguous subsequences of length $l$,
$$||P(S_t\otimes I)P(S_{t-1}\otimes I)\cdots  P(S_1\otimes I)P||\leq \lambda^{(t-l\omega)}$$
\label{nmeets2xxx}\end{theorem}

Compared with Theorem 3 of \cite{tac} in which $\omega = (m-1)^2$, it can be seen that
Theorem \ref{nmeets2xxx} leads to an improved lower bound on the convergence rate.

}

{\color{black}

\subsection{Least Squares Solutions} \label{least}

In the case when $Ax = b$ does not have a solution,
it is desirable to compute a least squares solution,
that is a value of $x$ for which $A'Ax = A'b$.
To obtain such a solution in a distributed manner, there are
two approaches to apply the idea of algorithm \rep{a1}.
One approach is to simply let each agent know a subset of the rows of the matrix
$\matt{A'A & A'b}$ instead of the matrix $\matt{A & b}$,
and then the algorithm for each agent needs to be modified accordingly;
apparently, this approach requires a centralized computation of
$A'A$ and $A'b$ before the algorithm is launched.
The other approach is described in Section IX of \cite{tac},
which needs a network wide design step and requires
an $mn$-dimensional state vector transmitted between neighboring agents at each
clock time.
Although both approaches have shortcomings, they are capable of finding
a least squares solution exponentially fast for a time-varying directed neighbor
graph sequence under appropriate joint connectedness.
The weakest possible conditions for exponential convergence
can be straightforwardly derived using the parameter-dependent notion of
graph connectivity introduced in this paper.

A least squares solution can be obtained using the existing distributed
convex optimization algorithms such as those in \cite{review12,review14,reviewlong,nedic,nedic2,nedic5}.
The algorithms in \cite{review12,review14,reviewlong} work only for time-invariant graphs.
Although the algorithms in \cite{nedic,nedic2,nedic5} work for time-varying graphs,
their convergence rates are slower than exponential due to the
vanishing step size.
The distributed least squares problem can also be solved by using
distributed averaging to compute the average of the matrix pairs
$(A'_iA_i, A'_ib_i)$. One drawback of this approach, however, is that the
amount of data to be communicated between agents does not
scale well with the number of agents.
See Section II in \cite{tac} for more comparisons of the existing approaches
to the distributed least squares problem.
Simply put, there is no distributed algorithm which is capable of finding a least squares solution exponentially fast
for time-varying directed graphs using at most an $n$-dimensional state
vector transmitted between neighboring agents at each clock time.
Whether or not such algorithms can be devised remains to be seen.

}

\section{Conclusions}

This paper has studied the stability of a distributed algorithm for solving linear algebraic equations.
A new notion of graph connectivity has been introduced, based on which
necessary and sufficient graph-theoretic conditions have been obtained for
the system defined by the algorithm to be exponentially stable under the most general assumption.
Such a parameter-dependent notion of graph connectivity has been proposed
for the first time here, and is expected to
have broader applications in distributed control problems
such as constrained consensus and distributed optimization \cite{nedic2}.

\section{Appendix} \label{appendix}

To prove Proposition \ref{ns}, we need the following lemma.

\begin{lemma}
Let $S$ be an $m\times m$ stochastic matrix with positive diagonal entries.
Then, $P(S\otimes I)P$ is a discrete-time stability matrix if and only if
for each strongly connected component $\scr{E} \subset \scr{V}$
in $\gamma(S)$ which does not have neighbors, there holds
$$\bigcap_{j\in\scr{E}} \scr P_j = \{0\}$$
\label{nec}\end{lemma}

To prove Lemma \ref{nec}, we will need the following ideas from \cite{Gallager} and \cite{seneta}.
Let $\mathbb{G}$ be a directed graph with vertex set $\scr{V}$.
We say that a vertex $j$ is {\em reachable} from $i$ if either $j=i$ or
 if there is a directed path from $i$ to $j$.
We say that a vertex $i$ is {\em essential}  if $i$ is reachable  from all vertices which are reachable from $i$.
It is known that every directed graph has at least one essential vertex (see Lemma 10 of \cite{pieee}).

To proceed, let us say that
vertices $i$ and $j$ are
 {\em mutually reachable} if each is reachable from the other.
 Mutual reachability is
   an equivalence  relation on
  $\scr{V}$, which partitions $\scr{V}$ into a disjoint union of a
   finite number of  equivalence classes. Note that if
 $i$ is an essential vertex of $\mathbb{G}$, then every vertex  in the equivalence class
  of $i$ is also essential.
  Thus, every directed graph possesses at least one
   mutually reachable equivalence class whose members are all essential.
Note that a mutually reachable equivalence class of $\mathbb{G}$ is a strongly connected component
in $\mathbb{G}$ which does not have observers,
and a strongly connected graph possesses exactly one
mutually reachable equivalence class.

We also need the mixed matrix norm defined in \S \ref{unique}.
Recall that the  mixed matrix norm of $Q\in\R^{mn\times mn}$ is
$||Q|| = |\langle Q\rangle |_{\infty}$
 where $\langle Q\rangle $ is the $m\times m$ matrix in $\R^{m\times m}$  whose $ij$th entry is $|Q_{ij}|$.
Let us note that for any $m\times m$ nonnegative matrix $M$, there holds
$$\displaystyle\left\langle P(M\otimes I)P \right\rangle \leq M$$
where  for any real  matrices $X$ and $Y$ of the same size,  $X \leq Y$  means
that $Y-X$ is a nonnegative matrix.
The above inequality can be further extended. For any positive integer $k$,
 $$\displaystyle\left\langle (P(M\otimes I)P)^k \right\rangle \leq M^k$$
Therefore,
\eq{\displaystyle\left\|(P(M\otimes I)P)^k\right\| \leq \displaystyle\left|M^k\right|_{\infty}\label{rao}}

\noindent
{\bf Proof of Lemma \ref{nec}:}
In view of the preceding discussion,
 $\gamma(S)$ has at
least one mutually reachable equivalence class $\scr{E}$
 whose members are all essential.
 If $\scr{E}=\scr{V}$, then $\gamma(S)$  is strongly
 connected.
If  $\bigcap_{i=1}^m \scr P_i = \{0\}$, then
$P(S\otimes I)P$ is a discrete-time stability matrix (see Proposition 2 of \cite{cdc13.1}).
If $\bigcap_{i=1}^m \scr P_i\neq \{0\}$, then it is straightforward to verify that
$P(S\otimes I)P\bar z = \bar z$
where $\bar z= {\rm stack}\{z,z,\ldots,z\}$ and
$z$ is any nonzero vector in $\bigcap_{i=1}^m \scr P_i$.
Thus,  $P(S\otimes I)P$ has an eigenvalue at $1$. Therefore, in the case when $\gamma(S)$
is strongly connected, $P(S\otimes I)P$ is a discrete-time stability matrix if and only if
$\bigcap_{i=1}^m \scr P_i = \{0\}$.

Now suppose  that  $\scr{E} = \{i_1,i_2,\ldots,i_k\}$  is a strictly proper subset of
   $\scr{V}$. Let
 $\pi$ be any permutation map for which $\pi(i_j) = j,\;j\in\{1,2,\ldots,k\}$, and
 let $Q$ be the corresponding
 permutation matrix.
Then, the transformation $S\rightarrow QSQ'$ block triangularizes $S$;
specifically, there are matrices $S_1$, $S_2$ and $S_3$ such that
 $$QSQ' = \matt{S_1 & S_2\cr 0 &S_3}$$
Set $\bar Q = Q\otimes I$. Note that $\bar Q$ is a permutation matrix and
that $\bar QP\bar Q'$ is a block diagonal, orthogonal projection matrix
where the $j$th diagonal block is $P_{\pi(i_j)}$, $j\in\mathbf{k}$.
Since  $QSQ'$ is block triangular, there are matrices $A$, $B$ and $C$ such that
 $$\bar QP(S\otimes I)P\bar Q' = \matt{A & B\cr 0 &C}$$
 Thus, the spectrum of $P(S\otimes I)P$ consists of the eigenvalues of $A$ and $C$.

Set $P_{\scr{E}} = {\rm diagonal}\{P_{i_1},P_{i_2},\ldots,P_{i_k}\}$.
Then, it is straightforward to verify that
$$\matt{A & B} = P_{\scr{E}}\matt{S_1\otimes I & S_2 \otimes I}P_{\scr{E}}$$
Note that $\matt{S_1 & S_2}$ is a nonnegative matrix
whose row sums all equal $1$.
Since $\scr{E}$ is a mutually reachable equivalence class, $\gamma(S_1)$ is strongly connected.
If $S_2 = 0$, then $S_1$ is a stochastic matrix. Since $\gamma(S_1)$ is strongly connected,
in view of the proceeding, $A=P_{\scr{E}}(S_1\otimes I)P_{\scr{E}}$ is a discrete-time stability matrix if and only if
$\bigcap_{j\in\scr{E}} \scr P_{j} = \{0\}$.
If $S_2$ does not equal the zero matrix, then $S_1$ is not a stochastic matrix and
$S_1\leq S_1+S_2$.
In this case, $\scr{E}$ is a strongly connected component of $\gamma(S_1)$
which has at least one neighbor.
Since $\gamma(S_1)$ is strongly connected with positive diagonal entries, $S_1$ is a primitive matrix.
Note that $S_1+S_2$ is a stochastic matrix.
By the Perron-Frobenius Theorem (see Theorem 1.1 of \cite{seneta}),
all the eigenvalues of $S_1$  are strictly less than $1$ in magnitude, from which it is possible to show that
$A=P_{\scr{E}}(S_1\otimes I)P_{\scr{E}}$ is a discrete-time stability matrix.
Toward\footnote{
We thank Anup Rao (School of Computer Science, Georgia Institute of Technology)
for pointing out a flaw in the original version of this proof and suggesting how to fix it.
}
this end,
for any square matrix $T$, let $\rho(T)$ denote the spectral radius of $T$.
From the Gelfand formula \cite{gelfand} and \rep{rao},
\begin{eqnarray*}
\rho(A) &=& \lim_{k\rightarrow\infty}\displaystyle\left\|\displaystyle\left(P_{\scr{E}}(S_1\otimes I)P_{\scr{E}}\right)^k\right\|^{\frac{1}{k}} \\
&\leq& \lim_{k\rightarrow\infty}\displaystyle\left|S_1^k\right|_{\infty}^{\frac{1}{k}}
=\rho(S_1)<1
\end{eqnarray*}
which implies that $A$ is a discrete-time stability matrix if the corresponding
strongly connected component $\scr{E}$ has at least one neighbor.

Next we turn to the spectrum of matrix $C$.
Let $\bar \scr{E} = \{i_{k+1},\ldots,i_m\}$ be the complement of $\scr{E}$ in $\scr{V}$.
Set $P_{\bar\scr{E}} = {\rm diagonal}\{P_{\pi(i_{k+1})},P_{\pi(i_{k+2})},\ldots,P_{\pi(i_{m})}\}$.
Then, it is straightforward to verify that
Note that $S_3$ is a stochastic matrix and
$\gamma(S_3)$ is the subgraph of $\gamma(S)$ induced by $\bar \scr{E}$.

We can repeat the preceding procedure for matrix $C$ so that the spectrum of $C$ consists of
the eigenvalues of a set of block matrices. Each of those block matrices can be written as
$$P_{\scr{U}} \displaystyle\left(S_{\scr{U}}\otimes I \right)P_{\scr{U}}$$
where each $\scr{U}$ is a strongly connected component of $\gamma(S)$. The preceding analysis
applies to each of these matrices. Therefore, we have proved the lemma.
\hfill$\qed$

\noindent
{\bf Proof of Proposition \ref{ns}:}
Suppose first that  $P(F\otimes I)P$ is a discrete-time stability matrix.
To prove the necessity, suppose that, to the contrary, there is a partially populated subset
$\scr{E} = \{i_1,i_2,\ldots,i_k\}$
which has no neighbors  in $\scr{V}\setminus\scr{E}$.
Let $\pi$ be any permutation map for which $\pi(i_j) = j$, $j\in\mathbf{k}$, and
 let $Q$ be the corresponding
 permutation matrix.
Then, the transformation $F\rightarrow QFQ'$ block triangularizes $F$.
Set $\bar Q = Q\otimes I$. Note that $\bar Q$ is a permutation matrix and
that $\bar QP\bar Q'$ is a block diagonal, orthogonal projection matrix
where the $j$th diagonal block is $P_{\pi(i_j)}$, $j\in\mathbf{k}$.
Since  $QFQ'$ is block triangular, there are matrices $A$, $B$ and $C$ such that
 $$\bar QP(F\otimes I)P\bar Q' = \matt{A & 0\cr B &C}$$
which implies that the spectrum of $P(F\otimes I)P$ consists of the eigenvalues of $A$ and $C$.
Let $F_{\scr{E}}$ be that $k\times k$ submatrix of $F$ whose $pq$th entry is
the $i_pi_q$th entry of $F$ for all $p,q\in\mathbf{k}$.
In other words, $F_{\scr{E}}$ is that submatrix of $F$ obtained by deleting rows and
columns whose indices are not in $\scr{E}$.
Since $F$ is a stochastic matrix and there are no arcs from $\scr{V}\setminus\scr{E}$ to $\scr{E}$,
it follows that $F_{\scr{E}}$ is a stochastic matrix and
its graph $\gamma(F_{\scr{E}})$ is the subgraph of $\gamma(F)$ induced by $\scr{E}$.
Set $P_{\scr{E}} = {\rm diagonal}\{P_{i_1},P_{i_2},\ldots,P_{i_k}\}$.
Then, it is straightforward to verify that
$$A = P_{\scr{E}} \displaystyle\left(F_{\scr{E}}\otimes I\right)P_{\scr{E}}$$
Since $P(F\otimes I)P$ is a discrete-time stability matrix, all the eigenvalues of $A$
are strictly less than $1$ in magnitude.
Since $\scr{E}$  is a partially populated subset, it follows by definition that
$\bigcap_{j\in\scr{E}} \scr P_j \neq \{0\}$. Then, for any nonzero vector $z$ in $\bigcap_{j\in\scr{E}} \scr P_j $, there holds
$A \bar z = \bar z$,
where $\bar z= {\rm stack} \{z,z,\ldots,z\}$.
Thus,  $A$ has an eigenvalue at $1$. But this
contradicts the fact that all the eigenvalues of $A$
are strictly less than $1$ in magnitude. Therefore, $\scr{E}$ must
have at least one neighbor in $\scr{V}\setminus\scr{E}$.

Now we turn to the proof of sufficiency.
Suppose that every partially populated $\scr{E}\subset\scr{V}$  has at least one neighbor in $\scr{V}\setminus\scr{E}$.
By Lemma \ref{nec}, $P(F\otimes I)P$ is a discrete-time stability matrix if and only if
for each strongly connected component $\scr{E}$ of $\gamma(F)$ which does not have any neighbors
in $\scr{V}\setminus\scr{E}$, the intersection $\bigcap_{i\in\scr{E}} \scr P_{i} $ is the zero subspace.
Thus, to prove the lemma, it is enough to show that if
every partially populated subset of agents $\scr{E}$  has at least one neighbor which is not in $\scr{E}$,
then for each strongly connected component $\scr{E}$ of $\gamma(F)$ which does not have any neighbors
in $\scr{V}\setminus\scr{E}$, the intersection $\bigcap_{i\in\scr{E}} \scr P_{i} $ is the zero subspace.
Suppose that, to the contrary, there is a strongly connected component
$\scr{U} = \{i_1,i_2,\ldots,i_k\}$ of $\gamma(F)$ such that $\scr{U}$  does not have neighbors
in $\scr{V}\setminus\scr{U}$ and
$\bigcap_{j\in\scr{U}} \scr P_{j} \neq \{0\}$. Then, $\scr{U}$ is a partially populated subset.
But this contradicts the hypothesis that every partially populated subset has at least one neighbor in its complement.
Therefore, such a strongly connected component $\scr{U}$ does not exist.
\hfill$\qed$

%

\bibliographystyle{plain}        
\bibliography{consensus}

\end{document}